\documentclass[10pt]{article}
\usepackage{pifont}
\usepackage{amsfonts}
\usepackage{amsmath}
\usepackage{mathrsfs}
\usepackage{mathrsfs, amscd,amssymb,amsthm,amsmath,bm,graphicx,psfrag,subfigure}

\makeatletter

\renewcommand{\@seccntformat}[1]{{\csname the#1\endcsname}{\normalsize .}\hspace{.5em}}
\makeatother

\def \[{\begin{equation}}
\def \]{\end{equation}}
\newtheorem{thm}{Theorem}[section]
\newtheorem{prop}{Proposition}

\newtheorem{lem}[thm]{Lemma}

\newenvironment{wst}
{\setlength{\leftmargini}{1.5\parindent}
 \begin{itemize}
 \setlength{\itemsep}{-1.1mm}}
{\end{itemize}}

\begin{document}
\setlength{\baselineskip}{15pt}
\begin{center}{\large \bf
Characterization of tricyclic graphs with exactly two $Q$-main eigenvalues\footnote{Financially supported by the National Natural Science
Foundation of China (Grant No. 11071096) and the Special Fund for Basic Scientific Research of Central Colleges (CCNU11A02015).}}

\vspace{2mm}

{\large Shuchao Li\footnote{E-mail: lscmath@mail.ccnu.edu.cn (S.C.
Li), 251028455@qq.com (X. Yang)}, Xue Yang}\vspace{2mm}

Faculty of Mathematics and Statistics,  Central China Normal
University, Wuhan 430079, P.R. China\vspace{1mm}

\end{center}
\vspace{4mm}

\noindent {\bf Abstract}: The signless Laplacian matrix of a graph $G$ is defined to be the sum of its adjacency
matrix and degree diagonal matrix, and its eigenvalues are called $Q$-eigenvalues of $G$. A
$Q$-eigenvalue of a graph $G$ is called a $Q$-main eigenvalue if it has an eigenvector the sum
of whose entries is not equal to zero. Chen and Huang [L. Chen, Q.X. Huang, Trees, unicyclic graphs and bicyclic graphs with exactly two $Q$-main eigenvalues, submitted for publication] characterized all trees, unicylic graphs and bicyclic graphs
with exactly two main $Q$-eigenvalues, respectively. As a continuance of it, in this paper, all tricyclic
graphs with exactly two $Q$-main eigenvalues are characterized.

\vspace{2mm} \noindent{\it Keywords}: Signless Laplacian; $Q$-Main eigenvalue; Tricyclic graph%Semi-edge walk matrix.

\vspace{2mm}


\vspace{4mm}

 {\setcounter{section}{0}
\section{\normalsize Introduction}\setcounter{equation}{0}
Let $G=(V_G,E_G)$ be a simple connected graph with vertex set
$V_G=\{v_1,\ldots,v_n\}$ and edge set $E_G\not=\emptyset$. The adjacency matrix
$A=A(G)=(a_{ij})$ of $G$ is an $n\times n$ symmetric matrix with $a_{ij}=1$ if and only if $v_i,
v_j$ are adjacent and 0 otherwise. Since $G$ has no loops, the main
diagonal of $A$ contains only 0's. Suppose the valence or
degree of vertex $v_i$ equals $d_G(v_i)$ (or $d_i$) for $i=1,\ldots,n$, and let
$D=D(G)$ be the diagonal matrix whose $(i,i)$-entry is $d_i,
i=1,2,\ldots, n$.  The matrix $Q(G)=D(G)+A(G)$ has been called the
\textit{signless Laplacian matrix} of $G$. Recently, this matrix attracts more and more researchers' attention.
For survey papers on this matrix the reader is referred to \cite{C-R-S2,C-R-S1,C-R-S3}. The
eigenvalues and $Q$-eigenvalues of $G$ are those of $A(G)$ and $Q(G)$,
respectively. An eigenvalue ($Q$-eigenvalue) of a graph $G$ is called a
\textit{main eigenvalue} ($Q$-\textit{main eigenvalue}) if it has an eigenvector the sum
of whose entries is not equal to zero. The Perron-Frobenius theorem
implies that the largest eigenvalue and $Q$-eigenvalue of $G$ are always
main.

A vertex of a graph $G$ is said to be \textit{pendant} if its degree is one. Let $PV(G)$ be the set of all pendants of $G.$
Throughout the text we denote by $P_n$ and $C_n$ the path and cycle on
$n$ vertices, respectively. Denote by $N(v)$ (or $N_G(v)$) the set of all
neighbors of $v$ in $G$. $c(G) = |E_G|-|V_G|+1$ is said to be \textit{cyclomatic number} of a connected
graph $G$. In particular, $G$ will be a tree, unicyclic graph, bicyclic graph, or tricyclic graph if $c(G) = 0, 1, 2$ or 3. Based on \cite{010,11,12,13,14,15,16,17} we know that a tricyclic graph $G$ contains at least 3 cycles and at most 7 cycles, furthermore, there
do not exist 5 cycles in $G$. Denote the set of tricyclic graphs on $n$ vertices by $\mathscr{T}_n$. Then let
$\mathscr{T}_n=\mathscr{T}_n^3\bigcup\mathscr{T}_n^4\bigcup\mathscr{T}_n^6\bigcup\mathscr{T}_n^7,$
where $\mathscr{T}_n^i$ denotes the set of tricyclic graphs in
$\mathscr{T}_n$ with exact $i$ cycles for $i=3,4,6,7.$

Let $G$ be a connected graph, and $\widetilde{G}$ be the subgraph of
$G$ which is obtained from $G$ by deleting its pendant vertex (if any)
continuously until there is no any pendant vertex left. Obviously,
$\widetilde{G}$ is a connected proper subgraph of $G$ if $G$ has
non-empty pendant vertex set $PV(G)$ and $\widetilde{G}=G$
otherwise. We call $\widetilde{G}$ the \textit{base} of $G$, and the
vertex in $V_{\widetilde{G}}$ the \textit{internal vertex} of $G$. If $G$
contains a cycle with $PV(G)\neq\emptyset$, the longest path
$P_G=v_0v_1\ldots v_k$ between $PV(G)$ and $V_{\widetilde{G}}$ (i.e., $v_0\in
PV(G), v_i\notin V_{\widetilde{G}}\, (i=1,\ldots,k-1)$ and $v_k\in
V_{\widetilde{G}}$) is called the \textit{longest pendant path}.
$T_1, T_2, \ldots, T_{15}$ are all the bases of tricyclic graphs depicted in Fig. 1. Throughout the context, we use
$\widetilde{G}=T_i$ to mean that $\widetilde{G}$ has the same cycle arrangements and the same labelled vertices
as that of $T_i$ for $i=1,2,\ldots, 15.$
Call a path $P=u_0u_1\ldots u_k\, (k\geq1)$ an \textit{internal path} of $G$ if
$d_G(u_0), d_G(u_k)\ge 3$ and $d_G(u_i)=2$ for $1\leq i\leq k-1$.
Obviously, there are two types of internal paths: $u_0\neq
u_k(k\geq1)$ and $u_0=u_k(k\geq3)$. For convenience, in what follows
we call the former \textit{internal path} and the latter \textit{internal cycle}.
\begin{figure}[h!]
\begin{center}
\psfrag{a}{$C_{r_1}$}\psfrag{b}{$C_{r_2}$}\psfrag{c}{$C_{r_3}$}
  \psfrag{d}{$u$}\psfrag{e}{$v$}\psfrag{1}{$v_1$}\psfrag{2}{$v_2$}
  \psfrag{f}{$u_1$}\psfrag{g}{$u_2$}\psfrag{3}{$v_3$}\psfrag{4}{$v_4$}
  \psfrag{h}{$v_1$}\psfrag{i}{$v_2$}\psfrag{j}{$v_3$}
  \psfrag{k}{$T_1$}\psfrag{l}{$T_2$}\psfrag{m}{$T_3$}
  \psfrag{n}{$T_4$}\psfrag{o}{$T_5$}\psfrag{p}{$T_6$}
  \psfrag{q}{$T_7$}\psfrag{r}{$T_8$}\psfrag{s}{$T_9$}
  \psfrag{t}{$T_{10}$}\psfrag{u}{$T_{11}$}\psfrag{v}{$T_{12}$}
  \psfrag{w}{$T_{13}$}\psfrag{x}{$T_{14}$}\psfrag{y}{$T_{15}$}
  \includegraphics[width=160mm]{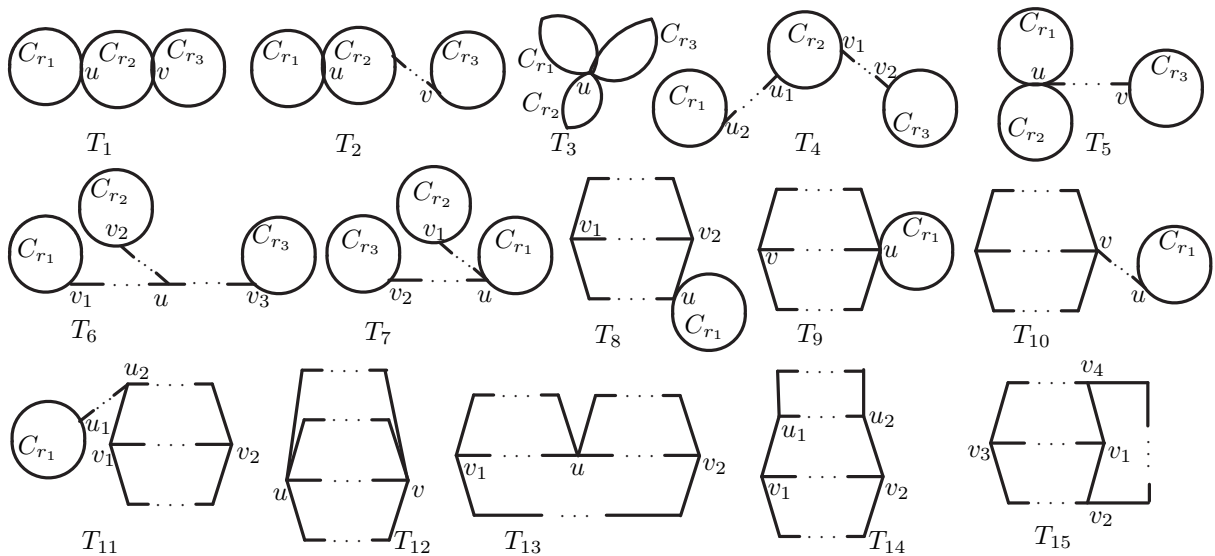}\\
  \caption{Graphs $T_1, T_2, \ldots, T_{15}$: all the bases of tricyclic graphs.}
\end{center}
\end{figure}

There are some literatures on main eigenvalues of adjacency matrix $A(G)$, there is a survey
in \cite{1} related to main eigenvalues of a graph. A long standing problem posed by Cvetkovi\'c (see \cite{2})
is to characterize graphs with exactly $k\, (k \ge 2)$ main eigenvalues. The graphs with $c(G) \le 3$
that have exactly two main eigenvalues are completely characterized (see, e.g. \cite{3,4,5,6,7,18}).
Motivated by these works, Chen and Huang \cite{10} obtained the following results on graphs with $Q$-main eigenvalues:
\begin{thm}[\cite{10}]
$G$ contains just one $Q$-main eigenvalue if and only if $G$ is regular.
\end{thm}
\begin{thm}[\cite{10}]
Let $G$ be a graph with signless Laplacian $Q$. Then $G$ has exactly two
$Q$-mian eigenvalues  if and only if there exist a unique pair of
integers a and b such that for any $v\in V_G$
\begin{equation}\label{eq:1.1}
\sum_{u\in N_G(v)}d_G(u)=ad_G(v)+b-d_G(v)^2
\end{equation}
\end{thm}

For convenience, let $\mathscr{G}_{a,b}$ be the set of the
connected graphs that satisfying (\ref{eq:1.1}), where $(a,b)$ is the
parameters of $G\in \mathscr{G}_{a,b}$. We know that $a>0$ and
$b\leq0$ from \cite{10}.
\begin{thm}[\cite{10}]
 Let $G\in\mathscr{G}_{a,b}$ with pendant vertex, then $a+b-1\geq2$.
 In particular, if $v\in V_G$ is a pendant vertex and u is the unique
 neighbor of v, then $d_G(u)=a+b-1.$
\end{thm}
\begin{thm}[\cite{10}]
Let  $G\in\mathscr{G}_{a,b}$ with pendant vertex, and
$P_{k+1}=v_0v_1v_2\ldots v_k$ be a longest pendant path in $G$. If $G$ contains
a cycle, then $1\leq k\leq2$. Moreover, if $k=2$, then $b\le -2$ and
$b^2a^2-(-2b^3+5b^2+b+d_{\widetilde{G}}(v_2))a+b^4-5b^3+5b^2+2b-d_{\widetilde{G}}(v_2)b+3d_{\widetilde{G}}(v_2)\leq0$.
\end{thm}

In particular, Chen and Huang \cite{10} characterized all trees, unicylic graphs and bicyclic graphs
with exactly two main $Q$-eigenvalues, respectively. As a continuance of it, in this paper we consider $Q$-main eigenvalues of tricyclic graphs.
We are to characterize all tricyclic graphs with exactly two $Q$-main eigenvalues in this paper. The organization of this work is
as follows: In Section 2, all tricyclic graphs without pandants having exactly two $Q$-main eigenvalues are determined.
In Section 3, all tricyclic graphs with pendants having exactly two $Q$-main eigenvalues are identified.

Further on we will need the following lemmas.

Based on (\ref{eq:1.1}) the next lemma follows immediately.
\begin{lem}\label{lem2}
If there exist $u,v\in V_G$ such that $d_G(u)=d_G(v)$, then
$\sum_{w\in N_G(u)}d_G(w)=\sum_{w\in N_G(v)}d_G(w)$. Moreover, if
$N_G(u)=\{u_1,u_2\}$ and $N_G(v)=\{v_1,v_2\}$, then
$d_G(u_1)+d_G(u_2)=d_G(v_1)+d_G(v_2).$
\end{lem}
\begin{lem}\label{lem1}
Let  $G\in\mathscr{G}_{a,b},$ $P=u_0u_1\ldots u_k$ be an internal
path or an internal cycle in $G$. Then
\begin{wst}
\item[{\rm (i)}] $k\leq3$, and

\item[{\rm (ii)}] If $k=3$, then there exists no internal path of length $2$, say $P_3=v_0v_1v_2$, in $G$
such that $d_G(v_0)=d_G(v_2)=d_G(u_0).$

\item[{\rm (iii)}] If $k=3$, then $d_G(u_0)=d_G(u_3)$. Moreover, if there exists another internal path $P_4'=v_0v_1v_2v_3$ in $G$,
then $d_G(v_0)=d_G(v_3)=d_G(u_0)=d_G(u_3).$
\end{wst}
\end{lem}
\begin{proof}
(i)\ \  On the contrary, suppose that $k\geq4$. By definition,
$d_G(u_0)=d_G(u_1)=d_G(u_2)=d_G(u_3)=2$. Note that
$d_{G}(u_1)=d_{G}(u_2)=2$, hence by Lemma \ref{lem2}, we have
$$
d_G(u_0)+2=\sum_{u\in N_G(u_1)}d_G(u)=\sum_{u\in
N_G(u_2)}d_G(u)=2+2=4,
$$
which gives $d_G(u_0)=2$, a contradiction to the condition that $d_G(u_0)\geq 3$.

(ii)\ \ If $k=3$, on the contrary, suppose that there exists an
internal path $P_3=v_0v_1v_2$ in $G$ such that
$d_G(v_0)=d_G(v_2)=d_G(u_0).$  Note that $d_{G}(u_1)=d_{G}(v_1)=2$,
hence by Lemma \ref{lem2}, we have
    $$
    d_G(u_0)+2=\sum_{u\in N_G(u_1)}d_G(u)=\sum_{u\in N_G(v_1)}d_G(u)=d_G(v_0)+d_G(v_2),
    $$
    which gives $d_G(u_0)=2$, a contradiction to the condition that $d_G(u_0)\geq 3$, as desired.

(iii)\ \  If $P$ is an internal cycle, it's trivial since $u_0=u_3$.
If $P$ is an internal path, note that $d_{G}(u_1)=d_{G}(u_2)=2$,
hence by Lemma \ref{lem2}, we have
    $$
    d_G(u_0)+2=\sum_{u\in N_G(u_1)}d_G(u)=\sum_{u\in N_G(u_2)}d_G(u)=d_G(u_3)+2,
    $$
    which gives that $d_G(u_0)=d_G(u_3)$. Moreover, if there exists another internal path $P_4'=v_0v_1v_2v_3$ in $G$, note that
$d_{G}(u_1)=d_{G}(v_1)=2$, hence by Lemma \ref{lem2}, we have
    $$
    d_G(u_0)+2=\sum_{u\in N_G(u_1)}d_G(u)=\sum_{u\in N_G(v_1)}d_G(u)=d_G(v_0)+2,
    $$
    which gives that $d_G(u_0)=d_G(v_0)$. Similarly, we can obtain that $d_G(u_3)=d_G(v_3)$. Then our result follows immediately.
\end{proof}

\begin{lem}\label{lem3}
Let $G\in\mathscr{G}_{a,b}$ be a graph containing a cycle and let $P_{k+1}=v_0v_1\ldots v_k$ be a longest pendant path in
$G$. Then we have $k=1$.
\end{lem}
\begin{proof}
By Theorem 1.4 we have $1\leq k\leq2$. Clearly, $d_G(v_1)=a+b-1$ by
Theorem 1.3. In order to complete the proof, it suffices to show
that $k\neq 2$. On the contrary, we assume that $k=2$, then $v_2\in
V_{\widetilde{G}}$. From Theorem 1.4, we have $b\le -2$ and
\begin{eqnarray}\label{eq:3.1}
&& b^2a^2-(-2b^3+5b^2+b+r)a+(b^4-5b^3+5b^2+2b-rb+3r)\leq 0,
\end{eqnarray}
where $r=d_{\widetilde{G}}(v_2)\geq2.$ The left sides of
(\ref{eq:3.1}) may be viewed as quadratic equation of $a$ and its
discriminant is
\begin{eqnarray*}
  \Delta &=& (-2b^3+5b^2+b+r)^2-4b^2(b^4-5b^3+5b^2+2b-rb+3r)\\
   &=& b^4+2b^3+b^2-2b^2r+2br+r^2.
\end{eqnarray*}

In order to complete the proof, it suffices to consider
$\Delta\geq0$. Note that $(b^2-b-r)^2-\Delta
=-4b^3>0,(b^2+b+r)^2-\Delta=4rb^2>0$, hence we have $
b^2-b-r>\sqrt{\Delta},\,b^2+b+r>\sqrt{\Delta}.$

By solving inequality (\ref{eq:3.1}), we have $a_1\leq a\leq a_2$,
where
\begin{eqnarray*}
  a_1 = \frac{-2b^3+5b^2+b+r-\sqrt{\Delta}}{2b^2}, \ \ \ \
  a_2= \frac{-2b^3+5b^2+b+r+\sqrt{\Delta}}{2b^2}.
\end{eqnarray*}
Notice that
\begin{eqnarray*}
 && a_1=\frac{-2b^3+5b^2+b+r-\sqrt{\Delta}}{2b^2} > \frac{-2b^3+5b^2+b+r-(b^2+b+r)}{2b^2}=-b+2, \\[3pt]
 && a_2=\frac{-2b^3+5b^2+b+r+\sqrt{\Delta}}{2b^2} < \frac{-2b^3+5b^2+b+r+(b^2-b-r)}{2b^2}=-b+3,
\end{eqnarray*}
hence $-b+2<a<-b+3$, we may get a contradiction since both of $a$ and $b$
are integers.

Hence, we obtain $k=1$, as desired.% It completes this proof.
\end{proof}

Given a graph $G\in\mathscr{G}_{a,b}$ with a cycle and pendants, if $u\in V_{\widetilde{G}}$ satisfies
$d_G(u)\neq d_{\widetilde{G}}(u)$, then $G$ must contain pendants attached to $u$ by Lemma \ref{lem3}.
Hence, $d_G(u)\in\{d_{\widetilde{G}}(u),a+b-1\}$ for any $u\in V_{\widetilde{G}}$, and $d_G(u)=1$ if $u\notin V_{\widetilde{G}}$.
For $u\in V_{\widetilde{G}}$, $u$ is called an \textit{attached vertex} if it joins some pendant vertices and \textit{non-attached vertex} otherwise.
Let $u\in V_G$ be an attached vertex, then $d_G(u)=a+b-1>d_{\widetilde{G}}(u)$. Applying (\ref{eq:1.1}) at $u$, we have $\sum_{v\in N_{\widetilde{G}}(u)}d_G(v)+(d_G(u)-d_{\widetilde{G}}(u))=ad_G(u)+b-d_G(u)^2$, which leads to $\sum_{v\in N_{\widetilde{G}}(u)}d_G(v)=-ab-b^2+2b+d_{\widetilde{G}}(u)$.
If $b=0$ then $\sum_{v\in N_{\widetilde{G}}(u)}d_G(v)=d_{\widetilde{G}}(u)$. On the other hand, since $d_G(v)\geq2$ for $v\in N_{\widetilde{G}}(u)$, we have $\sum_{v\in N_{\widetilde{G}}(u)}d_G(v)\geq2d_{\widetilde{G}}(u)$, a contradiction. Hence $b\le -1$. Thus the next lemma follows immediately.
\begin{lem}\label{lem4}
Let $G\in\mathscr{G}_{a,b}$ be a graph containing a cycle and pendants.
\begin{wst}
\item[{\rm (i)}]If $u\in
V_{\widetilde{G}}$, then $d_G(u)\in\{d_{\widetilde{G}}(u),a+b-1\}$ and $d_G(u)=1$ otherwise. %if $u\notin V_{\widetilde{G}}.$
\item[{\rm (ii)}]If $u\in V(\widetilde{G})$ is an attached vertex, then
$$\text{$b\le -1,\ \ \ d_G(u)=a+b-1>d_{\widetilde{G}}(u),\ \ \ \sum_{v\in
N_{\widetilde{G}}(u)}d_G(v)=-ab-b^2+2b+d_{\widetilde{G}}(u)$.}$$
\end{wst}
\end{lem}

\section{\normalsize Tricyclic graphs without pendants having exactly two $Q$-main eigenvalues}
In this section, we identify all the tricyclic graphs without pendants having exactly two $Q$-main eigenvalues.
\begin{figure}[h!]
\begin{center}
\psfrag{a}{$G_1\in\mathscr{G}_{8,-6}$}\psfrag{b}{$G_2\in\mathscr{G}_{7,-4}$}
\psfrag{c}{$G_3\in\mathscr{G}_{9,-6}$}
\psfrag{e}{$G_4\in\mathscr{G}_{7,-5}$}\psfrag{f}{$G_5\in\mathscr{G}_{6,-3}$}
\psfrag{g}{$G_6\in\mathscr{G}_{6,-3}$}\psfrag{h}{$G_7\in\mathscr{G}_{8,-6}$}
\psfrag{i}{$G_8\in\mathscr{G}_{7,-5}$}\psfrag{j}{$G_{9}\in\mathscr{G}_{6,-3}$}
\psfrag{l}{$G_{10}\in\mathscr{G}_{7,-2}$}\psfrag{m}{$G_{11}\in\mathscr{G}_{8,-6}$}
\psfrag{n}{$G_{12}\in\mathscr{G}_{6,0}$}\psfrag{o}{$G_{13}\in\mathscr{G}_{7,-4}$}
\psfrag{p}{$G_{14}\in\mathscr{G}_{7,-4}$}\psfrag{q}{$G_{15}\in\mathscr{G}_{6,-2}$}
\psfrag{r}{$G_{16}\in\mathscr{G}_{6,-2}$}\psfrag{s}{$G_{17}\in\mathscr{G}_{5,0}$}
\psfrag{t}{$G_{18}\in\mathscr{G}_{8,-7}$}\psfrag{u}{$G_{19}\in\mathscr{G}_{7,-5}$}
\psfrag{v}{$G_{20}\in\mathscr{G}_{7,-5}$}\psfrag{w}{$G_{21}\in\mathscr{G}_{6,-3}$}
\psfrag{x}{$G_{26}\in\mathscr{G}_{5,0}$}\psfrag{y}{$G_{27}\in\mathscr{G}_{6,-3}$}
\psfrag{z}{$G_{22}\in\mathscr{G}_{7,-4}$}\psfrag{1}{$G_{23}\in\mathscr{G}_{8,-7}$}
\psfrag{2}{$G_{24}\in\mathscr{G}_{6,-2}$}\psfrag{3}{$G_{25}\in\mathscr{G}_{7,-5}$}
\includegraphics[width=140mm]{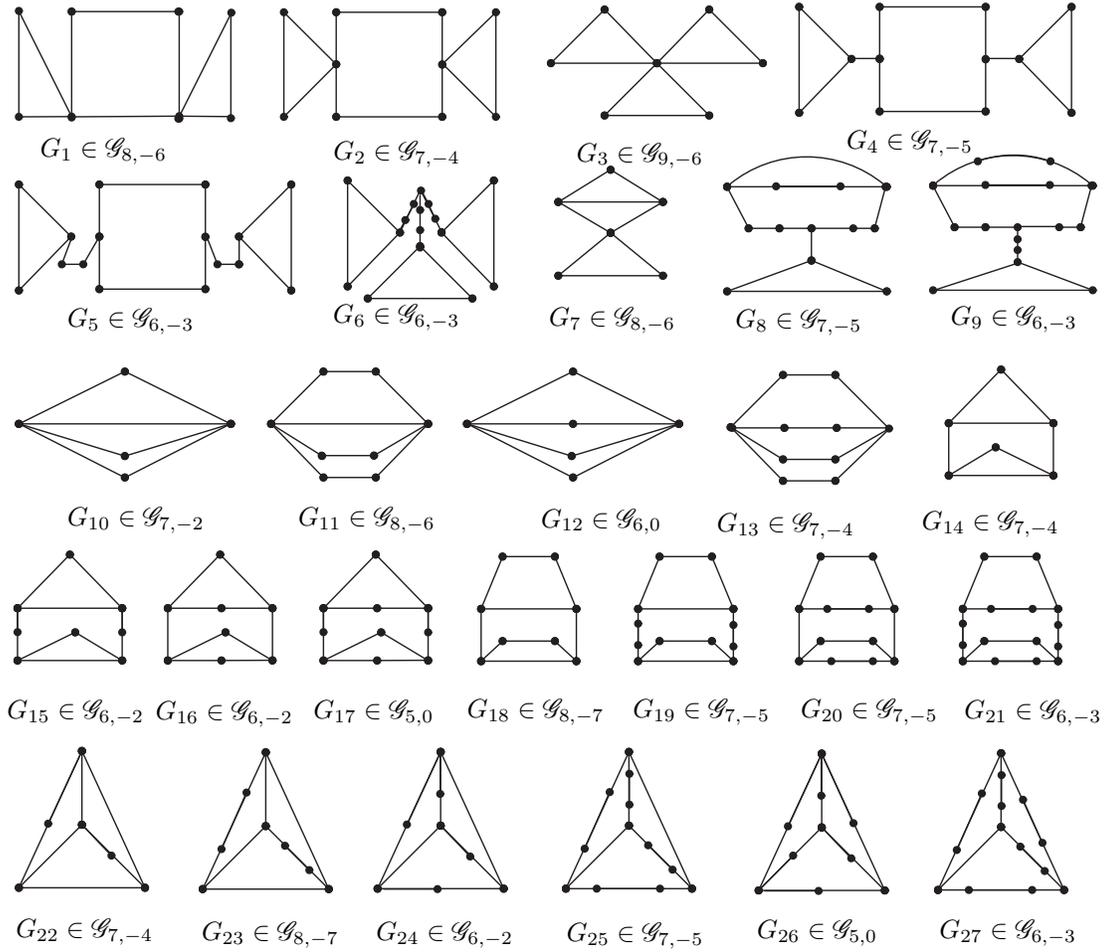}\\
\caption{Graphs $G_1, G_2, \ldots, G_{27}$.}
\end{center}
\end{figure}
\begin{thm}
$G_1\in\mathscr{G}_{8,-6}, G_2\in\mathscr{G}_{7,-4},
G_3\in\mathscr{G}_{9,-6}, G_4\in\mathscr{G}_{7,-5},
G_5\in\mathscr{G}_{6,-3}, G_6\in\mathscr{G}_{6,-3},
G_7\in\mathscr{G}_{8,-6}, G_8\in\mathscr{G}_{7,-5},
G_{9}\in\mathscr{G}_{6,-3}, G_{10}\in\mathscr{G}_{7,-2},
G_{11}\in\mathscr{G}_{8,-6}, G_{12}\in\mathscr{G}_{6,0},
G_{13}\in\mathscr{G}_{7,-4}, G_{14}\in\mathscr{G}_{7,-4},
G_{15}\in\mathscr{G}_{6,-2}, G_{16}\in\mathscr{G}_{6,-2},
G_{17}\in\mathscr{G}_{5,0}, G_{18}\in\mathscr{G}_{8,-7},
G_{19}\in\mathscr{G}_{7,-5}, G_{20}\in\mathscr{G}_{7,-5},
G_{21}\in\mathscr{G}_{6,-3}, G_{22}\in\mathscr{G}_{7,-4},
G_{23}\in\mathscr{G}_{8,-7}, G_{24}\in\mathscr{G}_{6,-2},
G_{25}\in\mathscr{G}_{7,-5}, G_{26}\in\mathscr{G}_{5,0},
G_{27}\in\mathscr{G}_{6,-3}$ (see Fig. 2) are all the tricyclic
graphs with no pendants having exactly two $Q$-main eigenvalues.
\end{thm}
\begin{proof}
Let $G\in \mathscr{G}_{a,b}$ be a tricyclic graph with no pendants having exactly two $Q$-main eigenvalues.
Hence, $G\cong \widetilde{G}\in\mathscr{T}_n=\mathscr{T}_n^3\bigcup\mathscr{T}_n^4\bigcup\mathscr{T}_n^6\bigcup\mathscr{T}_n^7$.
We consider the following possible cases.\vspace{2mm}

\noindent{\bf Case 1.}\ $G\in\mathscr{T}_n^3$. In this case, $G$
contains three cycles, say $C_{r_1}, C_{r_2}$ and $C_{r_3}$. By
Lemma \ref{lem1}(i), if $C_{r_i}$ is an internal cycle, then
$r_i=3$.

\medskip\noindent
$\bullet$ $G=T_1$; see Fig. 1. We have $r_1=r_3=3$. $C_{r_2}$
consists of the two internal paths: $P_{k_1+1}$ and $P_{k_2+1}$
connect $u,\,v$, where $d_G(u)=d_G(v)=4$. By Lemma \ref{lem1}(i), we
have $k_1\leq3,k_2\leq3$. By Lemma \ref{lem1}(ii), we have
$k_1\neq2,k_2\neq2$. Then without loss of generality, we assume
$k_1=1,k_2=3$ or $k_1=k_2=3$. It's simple to verify that $G\cong
G_1\in\mathscr{G}_{8,-6}$ if $k_1=1,k_2=3$, $G\cong
G_2\in\mathscr{G}_{7,-4}$ if $k_1=k_2=3$, where $G_1, G_2$ are
depicted in Fig. 2.

\medskip\noindent
$\bullet$ $G=T_2$; see Fig. 1. We have $r_1=r_3=3$. Suppose
$C_{r_1}=us_1s_2u$, and $C_{r_3}=vt_1t_2v$. Note that
$d_{G}(s_1)=d_{G}(t_1)=2$, hence by Lemma \ref{lem2}, we have
    $$
    2+4=d_{G}(s_2)+d_{G}(u)=d_{G}(t_2)+d_{G}(v)=2+3,
    $$
a contradiction.

\medskip\noindent
$\bullet$ $G=T_3$; see Fig. 1. We have $r_1=r_2=r_3=3$. It's routine to check that $G\cong G_3\in\mathscr{G}_{9,-6}$; see Fig. 2.

\medskip\noindent
$\bullet$ $G=T_4$; see Fig. 1. We have $r_1=r_3=3$. $C_{r_2}$
consists of the two internal paths: $P_{k_1+1}$ and $P_{k_2+1}$
connect $u_1$ and $v_1$; one internal path $P_{k_3+1}$ connects
$u_1, u_2$ and one internal path $P_{k_4+1}$ connects $v_1$ and
$v_2$, where $d_G(u_1)=d_G(u_2)=d_G(v_1)=d_G(v_2)=3$. By Lemma
\ref{lem1}(i), we have $k_i\leq3,\, i=1,2,3,4$. By Lemma
\ref{lem1}(ii), we have $k_i\neq2,\, i=1,2,3,4$. Then
$k_3=1,3,k_4=1,3$, furthermore, we assume, without loss of
generality, that $k_1=1,k_2=3$ or $k_1=k_2=3$. If $k_1=1,k_2=3$,
it's easy to prove that $k_3\neq1,\,3,k_4\neq1,\,3$. In fact, if
$k_3=1$, based on Lemma \ref{lem2}, we have
$$
2+2+3=\sum_{u\in N_G(u_2)}d_G(u)=\sum_{u\in N_G(u_1)}d_G(w)=2+3+3,
$$
a contradiction; if $k_3=3$,  based on Lemma \ref{lem2}, we have
$$
2+2+2=\sum_{u\in N_G(u_2)}d_G(u)=\sum_{u\in N_G(u_1)}d_G(w)=2+2+3,
$$
a contradiction. Hence, $k_3\neq1,\,3$. Similarly, $k_4\neq1,\,3$.
If $k_1=k_2=3$, it's impossible if $k_3=1,k_4=3$ or $k_3=3,k_4=1$.
In fact, if $k_3=1,k_4=3$, based on Lemma \ref{lem2}, we have
$$
2+2+3=\sum_{u\in N_G(u_1)}d_G(u)=\sum_{u\in N_G(v_1)}d_G(w)=2+2+2,
$$
a contradiction. Similarly, $k_3=3, k_4=1$ gives a contradiction. It
is routine to check that $G\cong G_4\in \mathscr{G}_{7,-5}$ if
$k_3=k_4=1$; whereas $G\cong G_5\in \mathscr{G}_{6,-3}$ if
$k_3=k_4=3.$

\medskip\noindent
$\bullet$ $G=T_5$; see Fig. 1. We have $r_1=r_2=r_3=3$. Suppose
$C_{r_1}=us_1s_2u$, and $C_{r_3}=vt_1t_2v$. Note that
$d_G(s_1)=d_G(t_1)=2$, based on Lemma \ref{lem2}, then we have
    $$
    2+5=d_{G}(s_2)+d_{G}(u)=d_{G}(t_2)+d_{G}(v)=2+3,
    $$
    which gives a contradiction.

\medskip\noindent
$\bullet$ $G=T_6$; see Fig. 1. We have $r_1=r_2=r_3=3$. There exist
three internal paths $P_{k_i+1}(i=1,2,3)$ connecting $u$ and $v_i$,
respectively. By Lemma \ref{lem1}(i), $k_1, k_2, k_3\leq3$. Since
$d_G(u)=d_G(v_1)=d_G(v_2)=d_G(v_3)=3$, by Lemma \ref{lem1}(ii),
$k_1, k_2, k_3\neq2$. Next we show that $k_1, k_2, k_3\neq1$. In
fact, suppose to the contrary, assume $k_1=1$. Based on Lemma
\ref{lem2}, we have
$$
2+2+3=\sum_{w\in N_G(v_1)}d_G(w)=\sum_{w\in
N_G(u)}d_G(w)=3+\sum_{w\in N_G(u)\backslash\{v_1\}}d_G(w),
$$
which is equivalent to $\sum_{w\in
N_G(u)\backslash\{v_1\}}d_G(w)=4$. It implies that $k_2=k_3=3$.
Based on Lemma \ref{lem2}, we have
$$
2+2+3=\sum_{w\in N_G(u)}d_G(w)=\sum_{w\in N_G(v_3)}d_G(w)=2+2+2,
$$
a contradiction. Hence, $k_1\not=1.$ Similarly, $k_2, k_3\not=1.$
Then $k_1=k_2=k_3=3$. It's simple to verify that $G\cong
G_6\in\mathscr{G}_{6,-3}$; see Fig. 2.

\medskip\noindent
$\bullet$ $G=T_7$; see Fig. 1. We have $r_1=r_2=r_3=3$. Suppose
$C_{r_1}=us_1s_2u$ and $C_{r_2}=v_1t_1t_2v_1$. Note that
$d_G(s_1)=d_G(t_1)=2$, based on Lemma \ref{lem2}, then we have
    $$
    2+4=d_{G}(s_2)+d_{G}(u)=d_{G}(t_2)+d_{G}(v)=2+3,
    $$
which gives a contradiction.\vspace{2mm}

\noindent{\bf Case 2.}\  $G\in\mathscr{T}_n^4$. In this case, it is easy to see that
$G$ contains an internal cycle $C_{r_1}$ with $r_1=3$.

\medskip\noindent
$\bullet$ $G=T_8$; see Fig. 1. Suppose $C_{r_1}=uz_1z_2u$. There
exist two internal paths $P_{k_1+1}$ and $P_{k_2+1}$ connecting
$v_1$ and $v_2$, one internal path $P_{k_3+1}$ connecting $u$ and
$v_1$, and one internal path $P_{k_4+1}$ connecting $u$ and $v_2$.
By Lemma \ref{lem1}(i), we have $k_i\leq3$ for $i=1, 2, 3, 4$. First we
show that $k_3=k_4=1$. By Lemma \ref{lem1}(iii), we have that $k_3\neq3$.
If $k_3=2$, note that $d_G(z_1)=d_G(x_1)=2$, then based on Lemma \ref{lem2},
we have
    $$
    2+4=d_{G}(z_2)+d_{G}(u)=d_{G}(v_1)+d_{G}(u)=3+4,
    $$
a contradiction. Thus, $k_3=1$; similarly, $k_4=1$. Note that
$d_G(u)=4,d_G(v_1)=d_G(v_2)=3$, hence by Lemma \ref{lem1}(iii), we
have that $k_1,k_2\neq3$. Without loss of generality, we assume
$k_1=1,k_2=2$ or $k_1=k_2=2$. It's simple to verify that $G\cong
G_7\in\mathscr{G}_{8,-6}$ if $k_1=1,k_2=2$. While if $k_1=k_2=2$,
applying (\ref{eq:1.1}) gives no integer solution, a contradiction.

\medskip\noindent
$\bullet$ $G=T_9$; see Fig. 1. There exist three internal paths
$P_{k_1+1}, P_{k_2+1}, P_{k_3+1}$ connecting $u$ and $v$. By Lemma \ref{lem1}(i), $k_1, k_2, k_3\leq3$. By Lemma \ref{lem1}(iii), $k_1,k_2,k_3\neq3$.
Without loss of generality we assume that $k_1=2$ and
$P_{k_1+1}=us_1v.$ Let $C_{r_1}=ut_1t_2u$. Note that
$d_G(t_1)=d_G(s_1)=2$, based on Lemma \ref{lem2}, we have
    $$
    2+5=d_{G}(t_2)+d_{G}(u)=d_{G}(v)+d_{G}(u)=3+5,
    $$
a contradiction.

\medskip\noindent
$\bullet$ $G=T_{10}$; see Fig. 1. By a similar discussion as in the
proof of $G=T_9$, we may obtain that there does not exist such graph
in $\mathscr{G}_{a,b},$ we omit the procedure here.

\medskip\noindent
$\bullet$ $G=T_{11}$; see Fig. 1. In this subcase, there exist two
internal paths $P_{k_1+1}$ and $P_{k_2+1}$ connecting $v_1,v_2$; one
internal path $P_{k_3+1}$ connecting $v_1,u_2$, one internal path
$P_{k_4+1}$ connecting $v_2,u_2$, and one internal path $P_{k_5+1}$
connecting $u_1,u_2$. Note that
$d_G(u_1)=d_G(u_2)=d_G(v_1)=d_G(v_2)=3$. By Lemma \ref{lem1}(i),
$k_i\leq3\, (i=1, 2, 3, 4, 5)$. By Lemma \ref{lem1}(ii),
$k_i\neq2\,(i=1, 2, 3, 4, 5)$. If $k_5=1$, then based on Lemma
\ref{lem2}, we have $k_3=k_4=3,k_1=1,k_2=3$ and it's simple to
verify that $G\cong G_8\in\mathscr{G}_{7,-5}$; see Fig. 2. If
$k_5=3$, similarly as above we obtain that $k_1=k_2=k_3=k_4=3$ and
it's  simple to verify that $G\cong G_{9}\in\mathscr{G}_{6,-3}$; see
Fig. 2.\vspace{2mm}

\noindent{\bf Case 3.}\ $G\in\mathscr{T}_n^6$. %In this case, $G\in \{T_{12}, T_{13}, T_{14}\}$ (see Fig. 1).

\medskip\noindent
$\bullet$ $G=T_{12}$; see Fig. 1. In this subcase, $G$ consists of
four internal paths: $P_{k_1+1},P_{k_2+1},P_{k_3+1},P_{k_4+1}$
connecting $u,v$. By Lemma \ref{lem1}(i), $k_i\leq3$ for $i=1, 2, 3,
4$. Note that $d_G(u)=d_G(v)=4$. By Lemma \ref{lem1}(ii), if there
exist $i_0\in\{1,2,3,4\}$ such that $k_{i_0}=3$, then $k_i\not=2$
for each $i\in \{1,2,3,4\}\setminus\{i_0\}$. Without loss of
generality, we may assume $k_1=1,k_2=k_3=k_4=2$, or
$k_1=1,k_2=k_3=k_4=3$, or $k_1=k_2=k_3=k_4=2$, or
$k_1=k_2=k_3=k_4=3$. It's simple to verify that $G\cong
G_{10}\in\mathscr{G}_{7,-2}$ if $k_1=1,k_2=k_3=k_4=2$; $G\cong
G_{11}\in\mathscr{G}_{8,-6}$ if $k_1=1, k_2=k_3=k_4=3$; $G\cong
G_{12}\in\mathscr{G}_{6,0}$ if $k_1=k_2=k_3=k_4=2$; $G\cong
G_{13}\in\mathscr{G}_{7,-4}$ if $k_1=k_2=k_3=k_4=3$, where $G_{10},
G_{11}, G_{12}$ and $G_{13}$ are depicted in Fig. 2.

\medskip\noindent
$\bullet$ $G=T_{13}$; see Fig. 1. By a similar discussion as in the
proof of $G=T_9$, we may obtain that there does not exist such graph
in $\mathscr{G}_{a,b},$ we omit the procedure here.

\medskip\noindent
$\bullet$ $G=T_{14}$; see Fig. 1. In this subcase, $G$ consists of
six internal paths: two paths $P_{k_1+1}, P_{k_2+1}$ connect
$u_1,u_2$, two paths $P_{k_3+1}, P_{k_4+1}$ connect $v_1,v_2$, one
path $P_{k_5+1}$ connects $u_1,v_1$ and one path $P_{k_6+1}$ connects
$u_2,v_2$. By Lemma \ref{lem1}(i), $k_i\leq3\, (i=1,2,3,4,5,6)$. Note that
$d_G(u_1)=d_G(u_2)=d_G(v_1)=d_G(v_2)=3$, by Lemma \ref{lem1}(ii), if there
exists $i_0\in\{1,2,3,4,5,6\}$ such that $k_{i_0}=3$, then
$k_i\neq2$ for each $i\in\{1,2,3,4,5,6\}\setminus \{i_0\}$.

If $k_1=1,k_2=2$, then $k_3=1,k_4=2$. Furthermore, if $k_5=1$, then
$k_6=1$, and it's simple to verify that $G\cong
G_{14}\in\mathscr{G}_{7,-4}$; If $k_5=2$, then $k_6=2$, and it's
simple to verify that $G\cong G_{15}\in\mathscr{G}_{6,-2}$, where
$G_{14}$ and $G_{15}$ are depicted in Fig. 2.

If $k_1=k_2=2$, then $k_3=k_4=2$. Furthermore, if $k_5=1$, then
$k_6=1$, and it's simple to verify that $G\cong
G_{16}\in\mathscr{G}_{6,-2}$; If $k_5=2$, then $k_6=2$, and it's
simple to verify that $G\cong G_{17}\in\mathscr{G}_{5,0}$, where
$G_{16}$ and $G_{17}$ are depicted in Fig. 2.

If $k_1=1,k_2=3$, then $k_3=1,k_4=3$. Furthermore, if $k_5=1$, then
$k_6=1$, and it's simple to verify that $G\cong
G_{18}\in\mathscr{G}_{8,-7}$; If $k_5=2$, then $k_6=2$, and it's
simple to verify that $G\cong G_{19}\in\mathscr{G}_{7,-5}$, where
$G_{18}$ and $G_{19}$ are depicted in Fig. 2.

If $k_1=k_2=3$, then $k_3=k_4=3$. Furthermore, if $k_5=1$, then
$k_6=1$, and it's simple to verify that $G\cong
G_{20}\in\mathscr{G}_{7,-5}$; If $k_5=3$, then $k_6=3$, and it's
simple to verify that $G\cong G_{21}\in\mathscr{G}_{6,-3}$, where
$G_{20}$ and $G_{21}$ are depicted in Fig. 2. \vspace{2mm}

\noindent{\bf Case 4.}\ $G\in\mathscr{T}_n^7$.\vspace{2mm}

In this case,  $G=T_{15}$ (see Fig. 1), hence $G$ consists of six
internal paths: one path $P_{k_{i-1}+1}$ connects $v_1,v_i\,
(i=2,3,4)$, one path $P_{k_4+1}$ connects $v_2,v_3$, one path
$P_{k_5+1}$ connects $v_3,v_4$, and one path $P_{k_6+1}$ connects
$v_2,v_4$. By Lemma \ref{lem1}(i), $k_i\leq3\, (i=1,2,3,4,5,6)$.  Note that
$d_G(v_1)=d_G(v_2)=d_G(v_3)=d_G(v_4)=3$. By Lemma \ref{lem1}(ii), if there
exists $i_0\in\{1,2,3,4,5,6\}$ such that $k_{i_0}=3$, then
$k_i\not=2$ for each $i\in\{1,2,3,4,5,6\}\backslash\{i_0\}$.

If $k_1=k_2=k_3=1$. It's simple to verify that $k_4=k_5=k_6=1$,
which implies that $G$ is a regular graph. By Theorem 1.1, $G$ contains
just one $Q$-main eigenvalue, a contradiction.

If $k_1=k_2=1$. Furthermore, if $k_3=2$, it's simple to verify that
$k_4=2,k_5=k_6=1$ and $G\cong G_{22}\in\mathscr{G}_{7,-4}$, if
$k_3=3$, it's simple to verify that $k_4=3,k_5=k_6=1$ and $G\cong
G_{23}\in\mathscr{G}_{8,-7}$. Here $G_{22}, G_{23}$ are depicted in
Fig. 2.

If $k_1=1$. Furthermore, if $k_2=k_3=2$, it's simple to verify that
$k_4=k_6=2,k_5=1$ and $G\cong G_{24}\in\mathscr{G}_{6,-2}$, if
$k_2=k_3=3$, it's simple to verify that $k_4=k_6=3,k_5=1$ and
$G\cong G_{25}\in\mathscr{G}_{7,-5}$. Here $G_{24}, G_{25}$ are
depicted in Fig. 2.

If $k_1=k_2=k_3=2$, it's simple to verify that $k_4=k_5=k_6=2$ and
$G\cong G_{26}\in\mathscr{G}_{5,0}$. If $k_1=k_2=k_3=3$, it's simple
to verify that $k_4=k_5=k_6=3$ and $G\cong
G_{27}\in\mathscr{G}_{6,-3}$. Here, $G_{26},G_{27}$ are depicted in
Fig. 2.

 This completes the proof.
\end{proof}

\section{\normalsize Tricyclic graphs with pendants having exactly two $Q$-main eigenvalues}
In this section, we identify all the tricyclic graphs with pendants
having exactly two $Q$-main eigenvalues.
%At first we will characterize the structure properties of tricyclic graphs with pendants having just two $Q$-main eigenvalues.
\setcounter{equation}{0}

\begin{lem}\label{lem5}
Given a tricyclic graph $G\in\mathscr{G}_{a,b}$ with pendants. If
$C_p=v u_1u_2 \ldots u_{p-1}v$ is an internal cycle of
$\widetilde{G}$ with $N_{\widetilde{G}}(v)=\{w,u_1, u_{p-1}\},$ then
$d_G(v)=3,\,p=3$ and $G\in\mathscr{G}_{6,-1}$. Moreover,
$d_G(u_1)=d_G(w)=2,\,d_G(u_{p-1})=a+b-1=4$.
\end{lem}
\begin{proof}
For convenience, let $r=d_G(w)$. Since $G$ is in $\mathscr{T}_n$, we
have $r\in \{2,3,4,5,a+b-1\}.$ Note that
$d_G(u_1),d_G(u_{p-1})\in\{2,a+b-1\}$, let
\[\label{3-01}
 t=|\{x: d_G(x)=2,\, x\in \{u_1, u_{p-1}\}\}|,
\]
hence, $t=0,1,2.$ We consider the following two possible cases
according to $d_G(v)$.\vspace{2mm}

\noindent\textbf{Case 1}. $d_G(v)=a+b-1>3$. In this case, by Lemma \ref{lem4} $b\le -1$, hence $a\geq 6$.\vspace{2mm}

If $d_G(u_1)=2$, applying (\ref{eq:1.1}) at $u_1$, we have
$d_G(u_2)=a-3\in\{2,a+b-1\}$. Notices that $a\geq6$, hence $a-3\neq
2$, i.e. $d_G(u_2)=a-3=a+b-1$, from which we get that $b=-2$; if
$d_G(u_1)=a+b-1$, applying (\ref{eq:1.1}) at $u_1 $, we have
$d_G(u_2)=-ab-b^2-a+b+3$.

Applying Lemma \ref{lem4}(ii) at $v$ yields
\begin{equation}\label{eq:3.2}
d_G(u_1)+d_G(u_{p-1})+d_G(w)=-ab-b^2+2b+3.
\end{equation}
In view of (\ref{3-01}) we have
\begin{equation}\label{eq:3.3}
2t+(a+b-1)(2-t)+r=-ab-b^2+2b+3,\ \ \ \ (t=0,1,2)
\end{equation}

$\bullet$ $t=0.$ In this subcase, (\ref{eq:3.3}) is equivalent to
$ab+b^2+2a=5-r$. As $d_G(u_1)=a+b-1$, we have
\begin{equation*}
\text{$ab+b^2+2a=5-r$  \ \  and \ \
    $d_G(u_2)=-ab-b^2-a+b+3$}.
\end{equation*}
with
\[\label{eq:3.4}
\text{$r\in
\{2,3,4,5,a+b-1\}$, \ \ \  $a+b>4$}
\]
Note that $d_G(u_2)\in\{2,a+b-1\}$. If $d_G(u_2)=2$, we get
$a+b=4-r$, a contradiction to (\ref{eq:3.4}); if $d_G(u_2)=a+b-1$,
we get $r=1$, a contradiction to (\ref{eq:3.4}) either.

$\bullet$ $t=1$. In this subcase, (\ref{eq:3.3}) gives
$ab+b^2+a-b=2-r$. Without loss of generality, we assume that
$d_G(u_1)=2$. Then we have $b=-2$ with $d_G(u_2)=a+b-1$. Hence,
$a=4+r$. Note that $a+b-1>3$, then we get $r\in \{3,4,5\}$. If
$r=3$, then $(a,b)=(7,-2)$ and $d_G(u_2)=4$. Applying (\ref{eq:1.1})
at $u_2$, it is easy to get that $d_G(u_3)=6\notin\{2,a+b-1=4\}$, a
contradiction. If $r=4$, then $(a,b)=(8,-2)$ and $d_G(u_2)=5$. We
apply (\ref{eq:1.1}) at $u_2$ to get that
$d_G(u_3)=8\notin\{2,a+b-1=5\}$, a contradiction. If $r=5$, then
$(a,b)=(9,-2)$ and $d_G(u_2)=6$. We apply (\ref{eq:1.1}) at $u_2$ to
get that $d_G(u_3)=10\notin\{2,a+b-1=6\}$, a contradiction.

$\bullet$ $t=2$. In this subcase, (\ref{eq:3.3}) gives
$ab+b^2-2b+1+r=0$. As $d_G(u_1)=2$, we have $b=-2$ with
$d_G(u_2)=a+b-1$. Hence, $2a=9+r$. Note that $a$ is an integer and
$a+b-1>3$, then we get $r=5$. Therefore, $(a,b)=(7,-2)$  and
$d_G(u_2)=4$. We apply (\ref{eq:1.1}) at $u_2$ to get that
$d_G(u_3)=6\notin\{2,a+b-1=4\}$, a contradiction. \vspace{2mm}

\noindent\textbf{Case 2}. $d_G(v)=3.$ In this case, $a+b-1\geq3$.
\vspace{2mm}

If $d_G(u_1)=2$, applying (\ref{eq:1.1}) at $u_1$, we get that
$d_G(u_2)=2a+b-7$; if $d_G(u_1)=a+b-1$, applying (\ref{eq:1.1}) at
$u_1$, we get that $d_G(u_2)=-ab-b^2+2b-1$.

Applying (\ref{eq:1.1}) at $v$ yields
\begin{equation}\label{eq:3.5}
d_G(u_1)+d_G(u_{p-1})+d_G(w)=3a+b-9.
\end{equation}
In view of (\ref{3-01}) we have
\begin{equation}\label{eq:3.6}
2t+(a+b-1)(2-t)+r=3a+b-9, \ \ \ \ (t=0, 1, 2)
\end{equation}

$\bullet$ $t=0$. By (\ref{eq:3.6}), $a-b=7+r$. Since
$d_G(u_1)=a+b-1$, we have
\begin{equation*}\label{eq:8}
\text{$ a-b=7+r$  \ \  and \ \
    $ d_G(u_2)=-ab-b^2+2b-1$}
\end{equation*}
with $d_G(u_2)\in\{2,a+b-1\}$. For $r\in\{2,a+b-1\}$, the equation
system gives that $(a,b)=(6,-3)$. However, $a+b-1=2<3$, a
contradiction. Moreover, it's easy to verify that the equation
system above has no integer solutions for $r\in\{3,4,5\}$, a
contradiction.

$\bullet$ $t=1$. By (\ref{eq:3.6}) $2a=10+r$. Note that $a$ is an
integer, then we have $r\in\{2,4,a+b-1\}$. Without loss of
generality, we assume that  $d_G(u_1)=2$, then
\begin{equation*}
\text{$ 2a=10+r$  \ \  and \ \
    $ d_G(u_2)=2a+b-7$}
\end{equation*}
with $d_G(u_2)\in\{2,a+b-1\}$. Since a and b are both integers
satisfying $a+b-1\geq3$, it's routine to check that only
$d_G(u_2)=a+b-1$ with $r=2$ holds, which gives that $(a,b)=(6,-1)$.
Then $d_G(u_2)=4$. Applying (\ref{eq:1.1}) at $u_2$, we have
$d_G(u_3)=3$. It implies that $u_3=v$ and $p=3$. Note that
$d_G(w)=r=2$.

$\bullet$ $t=2$. By (\ref{eq:3.6}) $3a+b=13+r$. Since $d_G(u_1)=2$,
we have
\begin{equation*}
\text{$ 3a+b=13+r$  \ \  and \ \
    $ d_G(u_2)=2a+b-7$}
\end{equation*}
with $d_G(u_2)\in\{2,a+b-1\}$.

For $r\in\{2,a+b-1\}$, it's routine to check that only
$d_G(u_2)=a+b-1$ with $r=a+b-1$ holds, which gives that $a=6$ and
$b\in\{-1,\,-2\}$ since a and b are both integers satisfying
$a+b-1\geq3$. First we consider $(a,b)=(6,-1)$. Then we have
$d_G(u_2)=4$. Applying (\ref{eq:1.1}) at $u_2$ yields
$d_G(u_3)=3\notin\{2,a+b-1=4\}$, a contradiction. Now we consider
$(a,b)=(6,-2)$. Then we have $d_G(u_2)=3$. Applying (\ref{eq:1.1})
at $u_2$ yields $d_G(u_3)=4\notin\{2,a+b-1=3\}$, a contradiction.

For $r=3$, if $d_G(u_2)=2$, then $(a,b)=(7,-5),a+b-1=1$, a
contradiction. So we have $d_G(u_2)=a+b-1$. Thus $(a,b)=(6,-2)$ and
$d_G(u_2)=3$. Applying (\ref{eq:1.1}) at $u_2$, we have
$d_G(u_3)=4\notin\{2,a+b-1=3\},$ a contradiction. For $r=4$, if
$d_G(u_2)=2$, then $(a,b)=(8,-7),a+b-1=0$, a contradiction. Hence,
we have $d_G(u_2)=a+b-1$,  thus $(a,b)=(6,-1)$ and $d_G(u_2)=4$.
Applying (\ref{eq:1.1}) at $u_2$, we have
$d_G(u_3)=3\notin\{2,a+b-1=4\}.$ For $r=5$, if $d_G(u_2)=2$, then
$(a,b)=(9,-9)$, a contradiction. Hence, we have $d_G(u_2)=a+b-1$,
thus $(a,b)=(6,0)$, a contradiction to $b\leq-1$.

This completes the proof.
\end{proof}
\begin{lem}
Let $G\in\mathscr{G}_{a,b}$ be a tricyclic graph with pendants. If
$C_p=v u_1 \ldots u_{p-1}v$ is an internal cycle of $\widetilde{G}$
with $N_{\widetilde{G}}(v)=\{u_1, u_{p-1}, w_1, w_2\}$, then
$d_G(v)=4,\,p=3$ and $G\in\mathscr{G}_{7,-1}$. Moreover,
$d_G(u_1)=d_G(w_1)=d_G(w_2)=2,\,d_G(u_{p-1})=a+b-1=5$.
\end{lem}
\begin{proof}
For convenience, let $d_1 =d_G(w_1),d_2=d_G(w_2)$. Note that $G\in
\mathscr{T}_n$, we have $d_1,d_2\in \{2,3,4,a+b-1\}.$ We consider
the following two possible cases according to $d_G(v)$.

\vspace{2mm}
\noindent{\bf Case 1}.  $d_G(v)=a+b-1>4.$ Note that $b\le -1$, hence in this case, $a\geq7$. \vspace{2mm}

If $d_G(u_1)=2$, we apply (\ref{eq:1.1}) at $u_1$ to get
$d_G(u_2)=a-3\in\{2,a+b-1\}$. As $a\geq 7$, we obtain $a-3\neq 2$.
Thus, $d_G(u_2)=a-3=a+b-1$, i.e., $b=-2$; If $d_G(u_1)=a+b-1$,
applying (\ref{eq:1.1}) at $u_1$, we get $d_G(u_2)=-ab-b^2-a+b+3$.

Applying Lemma \ref{lem4}(ii) at $v$ yields
\begin{equation}\label{eq:3.7}
d_G(u_1)+d_G(u_{p-1})+d_1+d_2=-ab-b^2+2b+4.
\end{equation}

We first consider $d_1\in \{2,a+b-1\}$. Let
\[\notag
 t_1=|\{x: d_G(x)=2,\, x\in \{u_1, u_{p-1}, w_1\}\}|,
\]
hence $t_1=0,1,2,3.$ Together with (\ref{eq:3.7}) we have
\[\label{eq:3.8}
2t_1+(a+b-1)(3-t_1)+d_2=-ab-b^2+2b+4,\ \ \ \ \  (t_1=0, 1, 2, 3).
\]

$\bullet$ $t_1=0$. By (\ref{eq:3.8}), $ab+b^2+3a+b=7-d_2$. As
$d_G(u_1)=a+b-1$, we have
 \begin{equation*}
\text{$ab+b^2+3a+b=7-d_2$  \ \  and \ \
    $ d_G(u_2)=-ab-b^2-a+b+3$}.
\end{equation*}
with $d_G(u_2)\in\{2,a+b-1\}$. It's routine to check that there is
no integer solution such that $a+b-1>4$ since
$d_2\in\{2,3,4,a+b-1\}$, a contradiction.
%Note that $a+b-1>4$ and $a, b$ are both integers. For $d_G(u_2)=2$, we have $2(a+b)=6-d_2$, a contradiction; for $d_G(u_2)=a+b-1$, we have $a+b=3+d_2$, hence $d_2\neq a+b-1$; on the other hand, for $d_G(u_2)=a+b-1$, the equation system gives $d_2b+2+4d_2+b=0$, i.e. $b=\frac{-2-4d_2}{d_2+1}$. Since $b$ is an integer, it's routine to check that $d_2\notin \{2,3,4\}$, a contradiction.

$\bullet$ $t_1=1$. By (\ref{eq:3.8}),  $ab+b^2+2a=4-d_2$. If
$d_G(u_1)=2$, then $b=-2$, which gives $d_2=0$, a contradiction; If
$d_G(u_1)=a+b-1$, then we have
 \begin{equation*}
\text{$ab+b^2+2a=4-d_2$  \ \  and \ \
    $ d_G(u_2)=-ab-b^2-a+b+3$}.
\end{equation*}
with $d_G(u_2)\in\{2,a+b-1\}$. It's routine to check that there is
no integer solution such that $a+b-1>4$ since
$d_2\in\{2,3,4,a+b-1\}$, a contradiction.
%For $d_G(u_2)=2$, the equation system gives $a+b=3-d_2$. As $a+b-1>4$, we get a contradiction; while for $d_G(u_2)=a+b-1$, the equation system gives $d_2=0$, a contradiction.

$\bullet$ $t_1=2$. By (\ref{eq:3.8}), $ab+b^2+a-b=1-d_2$. Without
loss of generality, we assume that $d_G(u_1)=2$, then we have $b=-2$
with $d_G(u_2)=a+b-1$. So $a=5+d_2$. It's obvious that $d_2\neq
a+b-1$. Moreover, since $a+b-1>4$, we have $a\geq8$, which implies
$d_2\geq3$. Then $d_2\in \{3,4\}.$ If $d_2=3$, then $(a,b)=(8,-2)$
and $d_G(u_2)=5$.  We apply (\ref{eq:1.1}) at $u_2$ to get
$d_G(u_3)=8\notin\{2,a+b-1=5\}$, a contradiction. If $d_2=4$, then
$(a,b)=(9,-2)$ and $d_G(u_2)=6$. We apply (\ref{eq:1.1}) at $u_2$ to
get $d_G(u_3)=10\notin\{2,a+b-1=6\}$, a contradiction.

$\bullet$ $t_1=3$. By (\ref{eq:3.8}), $ab+b^2-2b+d_2+2=0$. Since
$d_G(u_1)=2$, we have $b=-2$. So $2a=10+d_2$. Since $a+b-1>4$ and
$a$ is an integer, it's easy to verify that $d_2\notin
\{2,3,4,a+b-1\}$, a contradiction.

Hence, we conclude that $d_1\notin \{2,a+b-1\}$. Similarly,
$d_2\notin \{2,a+b-1\}$. Hence, $d_1,d_2\in\{3,4\}$. Moreover, by
Fig. 1, we know that $d_1=d_2=3$ and $\widetilde{G}\cong T_7,T_8$.
Let
$$t_2=|\{x: d_G(x)=2, x\in \{u_1, u_{p-1}\}\}|,$$
hence, $t_2= 0,1,2.$
Together with (\ref{eq:3.7}) we have
\begin{equation}\label{eq:3.9}
2t_2+(a+b-1)(2-t_2)+3+3=-ab-b^2+2b+4,\ \ \ \ \  (t_2=0, 1, 2)
\end{equation}

If $t_2=0$, then (\ref{eq:3.9}) gives $ab+b^2+2a=0$. Since $d_G(u_1)=a+b-1$, we
have $d_G(u_2)=-ab-b^2-a+b+3\in\{2,a+b-1\}$, which gives no integer solution such that $a+b-1>4$, a
contradiction.

If $t_2=1$, then (\ref{eq:3.9}) gives $ab+b^2+a-b+3=0$. Without loss
of generality, we assume that $d_G(u_1)=2$, then we have $b=-2$ with
$d_G(u_2)=a+b-1$, which gives that $(a,b)=(9,-2)$ and $d_G(u_2)=6$.
We apply (\ref{eq:1.1}) at $u_2$ to get
$d_G(u_3)=10\notin\{2,a+b-1=6\}$, a contradiction.

If $t_2=2$, then (\ref{eq:3.9}) gives $ab+b^2-2b+6=0$. Since
$d_G(u_1)=2$, we have $b=-2$, which gives that $(a,b)=(7,-2)$, a
contradiction to the assumption $a+b-1>4$.\vspace{2mm}

\noindent\textbf{Case 2}. $d_G(v)=4$. In this case, $a+b-1\geq3$.\vspace{2mm}

If $d_G(u_1)=2$, we apply (\ref{eq:1.1}) at $u_1$ to get
$d_G(u_2)=2a+b-8$; if $d_G(u_1)=a+b-1$, we apply (\ref{eq:1.1}) at
$u_1$ to get $d_G(u_2)=-ab-b^2+2b-2$.

Applying (\ref{eq:1.1}) at $v$, we have
\begin{equation}\label{eq:10}
d_G(u_1)+d_G(u_{p-1})+d_G(w_1)+d_G(w_2)=4a+b-16.
\end{equation}

We first consider $d_1,d_2\in \{2,a+b-1\}$. Let
\[\notag
 t_1'=|\{x: d_G(x)=2,\, x\in \{u_1, u_{p-1}, w_1, w_2\}\}|,
\]
hence $t_1'=0,1,2,3,4.$ Together with (\ref{eq:10}) we have
$$
2t_1'+(a+b-1)(4-t_1')=4a+b-16,\ \ \ \ \ (t_1'=0, 1, 2, 3, 4)
$$
If $t_1'=0,1$, let $d_G(u_1)=a+b-1$; if $t_1'=2$, then
$d_G(u_1)\in\{2,a+b-1\}$; if $t_1'=3,4$, let $d_G(u_1)=2$. It's
routine to check that only $d_G(u_2)=a+b-1$ with $t_1'=3$ holds,
which gives that $a=7$ and $b\in\{-1,\,-2,\,-3\}$ since a and b are
both integers satisfying $a+b-1\geq3$.

For $(a,b)=(7,-2)$, we have $d_G(u_2)=4$. We apply (\ref{eq:1.1}) at
$u_2$ to get $d_G(u_3)=6\notin\{2,a+b-1=4\}$, a contradiction. For
$(a,b)=(7,-3)$, we have $d_G(u_2)=3$. We apply (\ref{eq:1.1}) at
$u_2$ to get $d_G(u_3)=6\notin\{2,4,a+b-1=3\}$, a contradiction.

For $(a,b)=(7,-1)$, we have $d_G(u_2)=5$. We apply (\ref{eq:1.1}) at
$u_2$ to get $d_G(u_3)=4$. It implies that $u_3=v$ and $p=3$.
Together with $t_1'=3$, we have that $d_G(w_1)=d_G(w_2)=2$.

Now consider $d_1\in \{2,a+b-1\}$ and $d_2\in\{3,4\}$.
 Let
\[\notag
 t_2'=|\{x: d_G(x)=2,\, x\in \{u_1, u_{p-1}, w_1\}\}|,
\]
hence $t_2'=0,1,2,3.$ Together with (\ref{eq:10}) we have
\begin{equation}\label{eq:11}
2t_2'+(a+b-1)(3-t_2')+d_2=4a+b-16\ \ \ \ \ (t_2'=0, 1, 2, 3).
\end{equation}

$\bullet$ $t_2'=0$. By (\ref{eq:11}), $a-2b=13+d_2$. Since
$d_G(u_1)=a+b-1$, we have
$$
\text{$a-2b=13+d_2$ \ \ \ and\ \ \ $d_G(u_2)=-ab-b^2+2b-2$}
$$
with $d_G(u_2)\in\{2,a+b-1\}$, which, respectively, implies no
integer solution since $d_2\in\{3,4\}$, a contradiction.

$\bullet$ $t_2'=1$. By (\ref{eq:11}),  $2a-b=16+d_2$. First
consider $d_G(u_1)=a+b-1$, then we have
$$
\text{$2a-b=16+d_2$ \ \ \ and\ \ \ $d_G(u_2)=-ab-b^2+2b-2$}
$$
with $d_G(u_2)\in\{2,a+b-1\}$, which, respectively, implies no
integer solution since $d_2\in\{3,4\}$, a contradiction. Now
consider $d_G(u_1)=2$, then we have
$$
\text{$2a-b=16+d_2$ \ \ \ and\ \ \ $d_G(u_2)=2a+b-8$},
$$
with $d_G(u_2)\in\{2,a+b-1\}$, which, respectively, implies no
integer solution since $d_2\in\{3,4\}$, a contradiction.

$\bullet$ $t_2'=2$. By (\ref{eq:11}), $3a=19+d_2$. Since
$a$ is an integer, we have $d_2\notin \{3,4\}$, a contradiction.

$\bullet$ $t_2'=3$. By (\ref{eq:11}), $4a+b=22+d_2$.
Since $d_G(u_1)=2$, we have
$$
\text{$4a+b=22+d_2$ \ \ \ and\ \ \ $d_G(u_2)=2a+b-8$}
$$
with $d_G(u_2)\in\{2,a+b-1\}$. For $d_G(u_2)=2$, the equation system
implies no integer solution such that $a+b-1\geq3$, a contradiction;
for $d_G(u_2)=a+b-1$, the equation system implies that $a=7$, then
$a+b-1=4-d_2$. If $d_2=3$, then $(a,b)=(7,-3)$ and $d_G(u_2)=3$.
Applying (\ref{eq:1.1}) at $u_2$, we get that
$d_G(u_3)=6\notin\{2,a+b-1=3\}$, a contradiction. If $d_2=4$, then
$(a,b)=(7,-2)$ and $d_G(u_2)=4$. We apply (\ref{eq:1.1}) at $u_2$ to
get that $d_G(u_3)=6\notin\{2,a+b-1=4\}$, a contradiction.

Finally, we consider $d_1,d_2\in\{3,4\}$. Moreover, by Fig. 1, we
have $d_1=d_2=3$ and $\widetilde{G}\cong T_7, T_8$. Let
\[\notag
 t_2''=|\{x: d_G(x)=2,\, x\in \{u_1, u_{p-1}\}\}|,
\]
hence $t_2''=0,1,2.$  Combining with
(\ref{eq:10}) we have
\begin{equation}\label{eq:12}
2t_2''+(a+b-1)(2-t_2'')+3+3=4a+b-16 \ \ \ \ (t_2''=0,1,2).
\end{equation}

If $t_2''=0$, then (\ref{eq:12}) gives $2a-b=20$. Since
$d_G(u_1)=a+b-1$, we have $d_G(u_2)=-ab-b^2+2b-2\in\{2,a+b-1\}$,
which gives no integer solution satisfying $a+b-1\geq3$, a
contradiction.

If $t_2''=1$, then (\ref{eq:12}) gives $a=\frac{23}{3}$, a contradiction.

If $t_2''=2$, then (\ref{eq:12}) gives $4a+b=26$. Since
$d_G(u_1)=2$, we have $d_G(u_2)=2a+b-8\in\{2,a+b-1\}$. For
$d_G(u_2)=2$, there is no integer solution such that $a+b-1\geq3$, a
contradiction. So $d_G(u_2)=a+b-1$, which gives that $(a,b)=(7,-2)$
and $d_G(u_2)=4$. We apply (\ref{eq:1.1}) at $u_2$ to get
$d_G(u_3)=6\notin\{2,a+b-1=4\}$, a contradiction.

This completes the proof.
\end{proof}
\begin{figure}[h!]
\begin{center}
\psfrag{1}{$G_{28}\in\mathscr{G}_{6,-1}$}\psfrag{2}{$G_{29}\in\mathscr{G}_{6,-1}$}
\psfrag{3}{$G_{30}\in\mathscr{G}_{6,-1}$}
\psfrag{4}{$G_{31}\in\mathscr{G}_{6,-1}$}
\includegraphics[width=140mm]{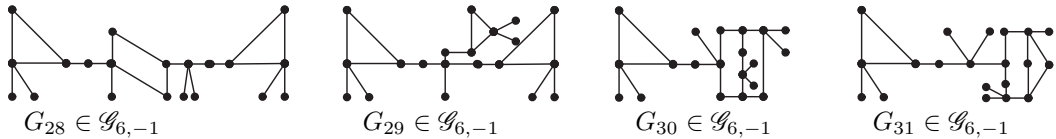}\\
\caption{Graphs $G_{28}, G_{29}, G_{30}$ and $G_{31}$.}
\end{center}
\end{figure}
\begin{prop}
Let $G\in\mathscr{G}_{a,b}$ be a tricyclic graph with pendants. If
$\widetilde{G}$ contains an internal cycle $C_p=w w_1 \ldots
w_{p-1}w$ with $d_{\widetilde{G}}(w)=3$ or $4$, then $G\cong
G_{28}\in\mathscr{G}_{6,-1},$ or  $G\cong
G_{29}\in\mathscr{G}_{6,-1},$ or  $G\cong
G_{30}\in\mathscr{G}_{6,-1},$ or  $G\cong
G_{31}\in\mathscr{G}_{6,-1}$; see Fig. 3.
\end{prop}
\begin{proof}
Based on Fig. 1, we obtain that $G\in\{T_1,\, T_2,\, T_4,\, T_5,\,
T_6,\, T_7,\, T_8,\, T_{10},\, T_{11}\}$. By Lemma \ref{lem5} and
Lemma 3.2, it's easy to see that $G\neq T_2,\,T_7$. \vspace{2mm}

\textbf{Case 1}. $G\in\mathscr{T}_n^3$. In this case,
$\widetilde{G}$ contains three cycles, say $C_{r_1}, C_{r_2}$ and
$C_{r_3}$. Then $G\in\{T_1,\, T_4,\, T_5,\, T_6\}$. By Lemma
\ref{lem5} and Lemma 3.2, if $C_{r_i}$ is an internal cycle of
$\widetilde{G}$, then $r_i=3$.

\medskip\noindent
$\bullet$ $G=T_1$; see Fig. 1. Then $r_1=r_3=3$. Note that
$d_{\widetilde{G}}(u)=d_{\widetilde{G}}(v)=4$. By Lemma 3.2, we
obtain that $G\in\mathscr{G}_{7,-1}$ and $d_G(u)=d_G(v)=4$. Denote
$C_{r_1}=ux_1x_2u$. Suppose $C_{r_2}$ consists of the two internal
paths of $\widetilde{G}$: $P_{k_1+1}=us_1s_2\ldots
s_{k_1}(s_{k_1}=v),\,P_{k_2+1}=ut_1t_2\ldots t_{k_2}(t_{k_2}=v)$
connects $u,\,v$. Then Lemma 3.2 implies that
$d_G(x_1)=d_G(s_1)=d_G(t_1)=2,d_G(x_2)=a+b-1=5$. Applying
(\ref{eq:1.1}) at $s_1$ and $t_1$ respectively yields
$d_G(s_2)=d_G(t_2)=5$. Lemma 3.2 implies that $d_G(v)=4$, hence,
$s_2,\,t_2\not=v$. We apply (\ref{eq:1.1}) at $s_2$ and $t_2$
respectively to get $d_G(s_3)=d_G(t_3)=4$, which implies that
$s_3=t_3=v$. However, by Lemma 3.2 we have $d_G(s_2)=d_G(t_2)=2$, a
contradiction.

\medskip\noindent
$\bullet$ $G=T_4$; see Fig. 1. Then $r_1=r_3=3$. Note that
$d_{\widetilde{G}}(u_2)=d_{\widetilde{G}}(v_2)=3$. By Lemma 3.1, we
obtain that $G\in\mathscr{G}_{6,-1}$ and $d_G(u)=d_G(v)=3$. Denote
$C_{r_1}=u_2x_1x_2u_2,\, C_{r_3}=v_2y_1y_2v_2$. Suppose the internal
paths of $\widetilde{G}$: $P_{k_1+1}=u_2q_1q_2\ldots
q_{k_1}(q_{k_1}=u_1)$ connects $u_2,\, u_1$,
$P_{k_2+1}=u_1s_1s_2\ldots s_{k_2}(s_{k_2}=v_1)$ and
$P_{k_3+1}=u_1t_1t_2\ldots t_{k_3}(t_{k_3}=v_1)$ connects
$u_1,\,v_1$, $P_{k_4+1}=v_1z_1z_2\ldots z_{k_4}(z_{k_4}=v_2)$
connects $v_1,\, v_2$. Then Lemma 3.1 implies that
$d_G(x_1)=d_G(y_1)=d_G(q_1)=2,\, d_G(x_2)=d_G(y_2)=a+b-1=4$.
Applying (\ref{eq:1.1}) at $q_1$ yields $d_G(q_2)=4$.

We first assume that $q_2=u_1$. We apply (\ref{eq:1.1}) at $u_1$ to
get $d_G(s_1)+d_G(t_1)=4$, which implies that $d_G(s_1)=d_G(t_1)=2$.
Applying (\ref{eq:1.1}) at $s_1$ and $t_1$ respectively yields
$d_G(s_2)=d_G(t_2)=3$. Hence, we have $s_2=t_2=v_1$ with
$d_G(v_1)=3$. Applying (\ref{eq:1.1}) at $v_1$, we get that
$d_G(z_1)=4$, moreover, $d_G(z_2)=2,\, d_G(z_3)=3$. Therefore,
$z_3=v$. Thus, it's easy to check that $G\cong
G_{28}\in\mathscr{G}_{6,-1}$; see Fig. 3.

Now we assume $q_2\not=u_1$. By the similar proof as in the case of
$q_2=u_1$, we may also get the graph $G\cong
G_{28}\in\mathscr{G}_{6,-1}.$

\medskip\noindent
$\bullet$ $G=T_5$; see Fig. 1. Then $r_3=3$. Note that
$d_{\widetilde{G}}(v)=3$. By Lemma 3.1, we obtain that
$G\in\mathscr{G}_{6,-1}$ and $d_G(v)=3$. Denote $C_{r_1}=vx_1x_2v$.
Suppose the internal paths of $\widetilde{G}$:
$P_{k_1+1}=vy_1y_2\ldots y_{k_1}(y_{k_1}=u)$ connects $v,\, u$. Then
Lemma 3.1 implies that $d_G(x_1)=d_G(y_1)=2,\, d_G(x_2)=a+b-1=4$.
Applying (\ref{eq:1.1}) at $y_1$ yields $d_G(y_2)=4$, moreover,
$d_G(y_3)=3\notin\{2,5,a+b-1=4\}$, a contradiction.

\medskip\noindent
$\bullet$ $G=T_6$; see Fig. 1. Then $r_1=r_2=r_3=3$. Note that
$d_{\widetilde{G}}(v_1)=d_{\widetilde{G}}(v_2)=d_{\widetilde{G}}(v_3)=3$.
By Lemma 3.1, we obtain that $G\in\mathscr{G}_{6,-1}$ and
$d_G(v_1)=d_G(v_2)=d_G(v_3)=3$. Denote $C_{r_1}=v_1x_1x_2v_1,\,
C_{r_2}=v_2y_1y_2v_2,\, C_{r_3}=v_3z_1z_2v_3$. Suppose the internal
paths of $\widetilde{G}$: $P_{k_1+1}=v_1s_1s_2\ldots
s_{k_1}(s_{k_1}=u)$ connects $v_1,\, u$, $P_{k_2+1}=ut_1t_2\ldots
t_{k_2}(t_{k_2}=v_2)$ connects $u,\, v_2$, $P_{k_3+1}=uq_1q_2\ldots
q_{k_3}(q_{k_3}=v_3)$ connects $u,\,v_3$. Then Lemma 3.1 implies
that $d_G(x_1)=d_G(y_1)=d_G(z_1)=d_G(s_1)=2,\,
d_G(x_2)=d_G(y_2)=d_G(z_2)=a+b-1=4$. Applying (\ref{eq:1.1}) at
$s_1$ yields $d_G(s_2)=4$.

We first assume that $s_2=u$. Then we have $d_G(u)=4$. We apply
(\ref{eq:1.1}) at $u$ to get that $d_G(t_1)+d_G(q_1)=4$, which
implies $d_G(t_1)=d_G(q_1)=2$. Applying (\ref{eq:1.1}) at $t_1$ and
$q_1$ respectively yields $d_G(t_2)=d_G(q_2)=3$. Hence, we have
$t_2=v_2,\, q_2=v_3$. Thus, it's easy to check that $G\cong
G_{29}\in\mathscr{G}_{6,-1}$; see Fig. 3.

Now we assume $s_2\not=u$. Applying (\ref{eq:1.1}) at $s_2$, we have
$d_G(s_3)=3$, which implies that $s_3=u$ with $d_G(u)=3$. Applying
(\ref{eq:1.1}) at $u$ yields $d_G(t_1)+d_G(q_1)=4$. Hence,
$d_G(t_1)=d_G(q_1)=2$. We apply (\ref{eq:1.1}) at $t_1$ to get that
$d_G(t_2)=4$, moreover, $d_G(t_3)=3$. Therefore, $t_3=v_2$. However,
Lemma 3.1 implies that $d_G(t_2)=2$, a contradiction. \vspace{2mm}

\textbf{Case 2}. $G\in\mathscr{T}_n^4$. In this case,
$\widetilde{G}$ contains an internal cycle. Then $G\in\{T_8,\,
T_{10},\, T_{11}$. By Lemma \ref{lem5} and Lemma 3.2, if $C_{r_i}$
is an internal cycle of $\widetilde{G}$, then $r_i=3$.

\medskip\noindent
$\bullet$ $G=T_8$; see Fig. 1. Then $r_1=3$. Note that
$d_{\widetilde{G}}(u)=4$. By Lemma 3.2, we obtain that
$G\in\mathscr{G}_{7,-1}$ and $d_G(u)=4$. Denote $C_{r_1}=ux_1x_2u$.
Suppose the internal paths of $\widetilde{G}$:
$P_{k_1+1}=us_1s_2\ldots s_{k_1}(s_{k_1}=v_1)$ connects $u,\, v_1$,
$P_{k_2+1}=ut_1t_2\ldots t_{k_2}(t_{k_2}=v_2)$ connects $u,\, v_2$,
$P_{k_3+1}=v_1y_1y_2\ldots y_{k_3}(y_{k_3}=v_2),\,
P_{k_4+1}=v_1z_1z_2\ldots z_{k_4}(z_{k_4}=v_2)$ connects
$v_1,\,v_2$. Then Lemma 3.2 implies that
$d_G(x_1)=d_G(s_1)=d_G(t_1)=2,\, d_G(x_2)=a+b-1=5$. We apply
(\ref{eq:1.1}) at $s_1$ and $t_1$ respectively to get
$d_G(s_2)=d_G(t_2)=5$. If $s_2\not=v_1$, we apply (\ref{eq:1.1}) at
$s_2$ to get $d_G(s_3)=4\notin\{2,3,a+b-1=5\}$, a contradiction.
Hence, $s_2=v_1$. Similarly, $t_2=v_2$. Applying (\ref{eq:1.1}) at
$v_1$ yields $d_G(y_1)+d_G(z_1)=5$, it's easy to see that this is
impossible.

\medskip\noindent
$\bullet$ $G=T_{10}$; see Fig. 1. Then $r_1=3$. Note that
$d_{\widetilde{G}}(u)=3$. By Lemma 3.1, we obtain that
$G\in\mathscr{G}_{6,-1}$ and $d_G(u)=3$. Denote $C_{r_1}=ux_1x_2u$.
Suppose the internal paths of $\widetilde{G}$:
$P_{k_1+1}=uy_1y_2\ldots y_{k_1}(y_{k_1}=v)$ connects $u,\, v$. Then
Lemma 3.1 implies that $d_G(x_1)=d_G(y_1)=2,\, d_G(x_2)=a+b-1=4$. We
apply (\ref{eq:1.1}) at $y_1$ to get $d_G(y_2)=4$. If $y_2\not=v$,
we apply (\ref{eq:1.1}) at $y_2$ to get
$d_G(y_3)=3\notin\{2,a+b-1=4\}$, a contradiction. Hence, $y_2=v$
with $d_G(v)=4$. Denote $N_G(v)=\{y_1,s_1,s_2,s_3\}$, then applying
(\ref{eq:1.1}) at $v$ yields $d_G(s_1)+d_G(s_2)+d_G(s_3)=5$, it's
easy to see that this is impossible.

\medskip\noindent
$\bullet$ $G=T_{11}$; see Fig. 1. Then $r_1=3$. Note that
$d_{\widetilde{G}}(u_1)=3$. By Lemma 3.1, we obtain that
$G\in\mathscr{G}_{6,-1}$ and $d_G(u_1)=3$. Denote
$C_{r_1}=u_1x_1x_2u_1$. Suppose the internal paths of
$\widetilde{G}$: $P_{k_1+1}=u_1y_1y_2\ldots y_{k_1}(y_{k_1}=u_2)$
connects $u_1,\, u_2$, $P_{k_2+1}=u_2s_1s_2\ldots
s_{k_2}(s_{k_2}=v_1)$ connects $u_2,\, v_1$,
$P_{k_3+1}=u_2t_1t_2\ldots t_{k_3}(t_{k_3}=v_2)$ connects $u_2,\,
v_2$, $P_{k_4+1}=v_1z_1z_2\ldots z_{k_4}(z_{k_4}=v_2),\,
P_{k_5+1}=v_1q_1q_2\ldots q_{k_5}(q_{k_5}=v_2)$ connects
$v_1,\,v_2$. Then Lemma 3.1 implies that $d_G(x_1)=d_G(y_1)=2,\,
d_G(x_2)=a+b-1=4$. Applying (\ref{eq:1.1}) at $y_1$, yields
$d_G(y_2)=4$.

We first consider $y_2=u_2$ with $d_G(u_2)=4$. Applying
(\ref{eq:1.1}) at $u_2$, we have $d_G(s_1)+d_G(t_1)=4$, which
implies that $d_G(s_1)=d_G(t_1)=2$. We apply (\ref{eq:1.1}) at $s_1$
and $t_1$ respectively to get $d_G(s_2)=d_G(t_2)=3$. It implies that
$s_2=v_1,\, t_2=v_2$ with $d_G(v_1)=d_G(v_2)=3$. Applying
(\ref{eq:1.1}) at $v_1$, yields $d_G(z_1)+d_G(q_1)=6$. Hence,
without loss of generality, we assume $d_G(z_1)=2,\, d_G(q_1)=4$.
Therefore, $d_G(z_2)=4,\, d_G(z_3)=3; d_G(q_2)=2,\, d_G(q_3)=3$.
Thus, $z_3=q_3=v_2$. It's easy to check that $G\cong
G_{30}\in\mathscr{G}_{6,-1}$; see Fig. 3.

Now we consider $y_2\not=u_2$. Then we apply (\ref{eq:1.1}) at $y_2$
to get $d_G(y_3)=3$, which implies that $y_3=u_2$ with $d_G(u_2)=3$.
Applying (\ref{eq:1.1}) at $u_2$ yields $d_G(s_1)+d_G(t_1)=4$.
Hence, $d_G(s_1)=d_G(t_1)=2$. Applying (\ref{eq:1.1}) at $s_1$ and
$t_1$ respectively, we have $d_G(s_2)=d_G(t_2)=4$. For the subcase
of $s_2=v_1$ with $d_G(v_1)=4$, we may get $d_G(z_1)+d_G(q_1)=4$ by
applying (\ref{eq:1.1}) at $v_1$. Hence, we have
$d_G(z_1)=d_G(q_1)=2$. We apply (\ref{eq:1.1}) at $z_1$ and $q_1$
respectively to get $d_G(z_2)=d_G(q_2)=3$. Therefore, $z_2=q_2=v_2$
with $d_G(v_2)=3$. It implies that $t_2\not=v_2$. Then applying
(\ref{eq:1.1}) at $t_2$ yields $d_G(t_3)=3$. Thus, $t_3=v_2$. It's
easy to check that $G\cong G_{31}\in\mathscr{G}_{6,-1}$; see Fig. 3.  For the
subcase of $s_2\not=v_1$, we may also get the graph  $G\cong
G_{31}\in\mathscr{G}_{6,-1}$ by the similar proof as
above.\vspace{2mm}

Thus, we complete the proof.
\end{proof}

\begin{prop}
Let $G$ be a tricyclic graph with pendants satisfying
$\widetilde{G}=T_3;$ see Fig. 1. Then $G\notin\mathscr{G}_{a,b}$.
\end{prop}
\begin{proof}For the graph $G$, let $u$ be the common vertex of $C_{r_1}, C_{r_2}$ and $C_{r_3}$ as depicted in $T_3$ of Fig. 1.
Assume that $N_{\widetilde{G}}(u)\cap V_{C_{r_1}}=\{u_1, u_{p-1}\},
N_{\widetilde{G}}(u)\cap V_{C_{r_2}}=\{w_1, w_2\},
N_{\widetilde{G}}(u)\cap V_{C_{r_3}}=\{w_3, w_4\}$. Suppose to the
contrary that $G\in\mathscr{G}_{a,b}$. Note that
$d_G(u_1), d_G(u_{p-1}), d_G(w_1), d_G(w_2), d_G(w_3), d_G(w_4)\in
\{2,a+b-1\}.$ Then we may let
\[\label{eq:3.24}
t=|x:d_G(x)=2, x\in \{u_1, u_{p-1}, w_1, w_2, w_3, w_4\}|.
\]
Hence, $0\leq t\leq 6$. Note that if $t=0$, we have $d_G(u_1)=a+b-1$;
if $t\in\{1,2,3,4,5,6\}$, without loss of generality, we can assume
$d_G(u_1)=2$. We proceed by distinguish the following two possible cases.
Note that
\[\label{eq:3.25}
 d_G(u_2)\in\{2,\ \ a+b-1\}.
\]

\vspace{2mm}

\textbf{Case 1}. $d_G(u)=a+b-1>6$. In this case as $b\le -1$, we have $a\geq 9$.\vspace{2mm}

If $d_G(u_1)=2$, applying (\ref{eq:1.1}) at $u_1$, we have
$d_G(u_2)=a-3\in\{2,a+b-1\}$. As $a\geq 9$, we get that $a-3\neq 2$.
Thus, $d_G(u_2)=a-3=a+b-1$, i.e., $b=-2$; if $d_G(u_1)=a+b-1$,
applying (\ref{eq:1.1}) at $u_1$, we have $d_G(u_2)=-ab-b^2-a+b+3$.
Apply Lemma \ref{lem4}(ii) at $u$ yields
$$
d_G(u_1)+d_G(u_{p-1})+d_G(w_1)+d_G(w_2)+d_G(w_3)+d_G(w_4)=-ab-b^2+2b+6.
$$
Together with (\ref{eq:3.24}) we have
$2t+(a+b-1)(6-t)=-ab-b^2+2b+6$. For $d_G(u_1)=2$, we have
\begin{equation*}\label{eq:3.26}
\left\{\begin{array}{ll}
  2t+(a+b-1)(6-t)=-ab-b^2+2b+6\ \ \ \ \ (t=1,2,3,4,5,6)\\
  b=-2
   \end{array}
 \right.
\end{equation*}
For $d_G(u_1)=a+b-1$, we have
\begin{equation*}\label{eq:3.27}
 \left\{\begin{array}{ll}
  2t+(a+b-1)(6-t)=-ab-b^2+2b+6\ \ \ \ \ (t=0)\\
  d_G(u_2)=-ab-b^2-a+b+3
   \end{array}
 \right.
\end{equation*}
It's routine to check that there is no integer solution satisfying
$a+b-1>6$, a contradiction.
%By (\ref{eq:3.26}), we get that $(t-4)(a-5)=4$. Together with  $a\ge
%9$, we have  $a-5=4,t-4=1$, thus $(a,b)=(9,-2),a+b-1=6$, a
%contradiction to the assumption $a+b-1>6$. On the other hand, in
%view of (\ref{eq:3.25}) and (\ref{eq:3.27}) with $t=0$, we get no
%integer solution satisfying $a+b-1\geq3$, a contradiction.
\vspace{2mm}

\textbf{Case 2}.  $d_G(u)=6$. In this case, $a+b-1\geq3$\vspace{2mm}

If $d_G(u_1)=2$, applying (\ref{eq:1.1}) at $u_1$, we have
$d_G(u_2)=2a+b-10$; if $d_G(u_1)=a+b-1$, applying (\ref{eq:1.1}) at
$u_1$, we have $d_G(u_2)=-ab-b^2+2b-4$.

Applying (\ref{eq:1.1}) at $u$ yields
$$
d_G(u_1)+d_G(u_{p-1})+d_G(w_1)+d_G(w_2)+d_G(w_3)+d_G(w_4)=6a+b-36.
$$
Together with (3.24) we have $2t+(a+b-1)(6-t)=6a+b-36$. For
$d_G(u_1)=2$, we have
\begin{equation}\label{eq:3.29}
\left\{\begin{array}{ll}
  2t+(a+b-1)(6-t)=6a+b-36\ \ \ \ \ (t=1,2,3,4,5,6)\\[3pt]
  d_G(u_2)=2a+b-10
   \end{array}
 \right.
\end{equation}
For $d_G(u_1)=a+b-1$, we have
\begin{equation}\label{eq:3.28}
 \left\{\begin{array}{ll}
   2t+(a+b-1)(6-t)=6a+b-36\ \ \ \ \ (t=0)\\[3pt]
  d_G(u_2)=-ab-b^2+2b-4
   \end{array}
 \right.
\end{equation}

In view of (\ref{eq:3.25}) and (\ref{eq:3.29}), if $d_G(u_2)=2$,
then we obtain that $(a-9)(10-t)=0$. It implies that $a=9$. Then
$b=-6$, hence $a+b-1=2<3$, a contradiction; if $d_G(u_2)=a+b-1$,
then $a=9$, so we get that $(5-t)(b+6)=0$, which implies that $t=5$
as $b\not=-6$. Thus $b\in\{-1,-2,-3,-4,-5\}$ since $b\le -1$ and
$a+b-1\geq3$.

Note that $t=5$, without loss of generality, we assume
$d_G(u_1)=d_G(u_{p-1})=2$. If $b=-1$, then $(a,b)=(9,-1)$ and
$d_G(u_2)=7$. Applying (\ref{eq:1.1}) at $u_2$, it's routine to
verify that $d_G(u_3)=6\notin\{d_{\widetilde{G}}(u_3)=2,a+b-1=7\}$,
a contradiction; for $b\in \{-2,-3,-4,-5\}$, we will also get a
contradiction by the similar proof as above.
%If $b=-2$, then $(a,b)=(9,-2)$, since
%$d_G(u_1)=2,d_G(u_2)=6$. Applying (\ref{eq:1.1}) at $u_2$, it's
%routine to verify that
%$d_G(u_3)=10\notin\{d_{\widetilde{G}}(u_3)=2,a+b-1=6\}$, a
%contradiction; If $b=-3$, as $d_G(u_1)=2,d_G(u_2)=a+b-1=5$, then
%$(a,b)=(9,-3)$. Applying (\ref{eq:1.1}) at $u_2$, it's routine to
%verify that $d_G(u_3)=12\notin\{d_{\widetilde{G}}(u_3)=2,a+b-1=5\}$,
%a contradiction; If $b=-4$, note that $d_G(u_1)=2,d_G(u_2)=a+b-1=4$,
%hence $(a,b)=(9,-4)$.  Applying (\ref{eq:1.1}) at $u_2$, it's
%routine to verify that
%$d_G(u_3)=12\notin\{d_{\widetilde{G}}(u_3)=2,a+b-1=4\}$, a
%contradiction; If $b=-5$, since $d_G(u_1)=2,d_G(u_2)=a+b-1=3$, we
%get $(a,b)=(9,-5)$. Applying (\ref{eq:1.1}) at $u_2$, it's routine
%to verify that
%$d_G(u_3)=10\notin\{d_{\widetilde{G}}(u_3)=2,a+b-1=3\}$,
% a contradiction.

From (\ref{eq:3.28}), we get that $(a,b)=(9,-6)$, hence $a+b-1=2<3$,
a contradiction.
\end{proof}
\begin{prop}
Let $G$ be a tricyclic graph with pendants. If $\widetilde{G}=T_9$ (see Fig. 1), then $G\notin\mathscr{G}_{a,b}$.
\end{prop}
\begin{proof}
Note that $\widetilde{G}=T_9$, hence $\widetilde{G}$ contains an
internal cycle $C_{r_1}:uu_1\ldots u_{p-1}u$ and three paths
$P_{k_1+1}:ux_1\ldots x_{k_1}(x_{k_1}=v), P_{k_2+1}:uy_1\ldots
y_{k_2}(y_{k_2}=v), P_{k_3+1}:uz_1\ldots z_{k_3}(z_{k_3}=v)$
connecting $u,v$; see $T_9$ of Fig.~1. Without loss of generality,
we assume $k_2\geq2,k_3\geq2$. Thus,
$d_G(u_1),d_G(u_{p-1}),d_G(y_1),d_G(z_1)\in\{2,a+b-1\}$. Denote
$N_{\widetilde{G}}(u)=\{u_1,u_{p-1},w,y_1,z_1\}$, where $w=v$ if
$k_1=1$; $w=x_1\not=v$ if $k_1\geq2$. Thus, we let
$r=d_G(w)\in\{2,3,a+b-1\}$ for convenience.  Let
\[\label{3-03}
t=|\{x: d_G(x)=2, x\in \{u_1, u_{p-1}, y_1, z_1\}\  \}|,
\]
hence, $0\le t\le 4$. Note that if $t=0,1$, let $d_G(u_1)=a+b-1$; if
$t=2$, we have $d_G(u_1)\in\{2,a+b-1\}$; if $t\in\{3,4\}$, let
$d_G(u_1)=2$.

In what follows, according to the vertex $u$ of maximum degree in
$\widetilde{G},$ we distinguish the following cases to prove our
result.

 \vspace{2mm}
 \textbf{Case 1}. $u$ is an attached vertex,
i.e. $d_G(u)=a+b-1>5$. In this case, as $b\le -1,$ we have
$a\geq8$.\vspace{2mm}

If $d_G(u_1)=2$, we apply (\ref{eq:1.1}) at $u_1$ to get
$d_G(u_2)=a-3\in\{2,a+b-1\}$. Since $a\geq8$, we get that $a-3\neq
2$. Hence, $d_G(u_2)=a-3=a+b-1$, i.e., $b=-2$; if $d_G(u_1)=a+b-1$,
we apply (\ref{eq:1.1}) at $u_1$ to get $d_G(u_2)=-ab-b^2-a+b+3$.

Applying Lemma \ref{lem4}(ii) at $u$, we have
\begin{equation}\label{eq:3.14}
\sum_{w\in N_{\widetilde{G}}(u)}=-ab-b^2+2b+5.
\end{equation}

Together with (\ref{3-03}), we have $r+2t+(a+b-1)(4-t)=-ab-b^2+2b+5.
$ For $d_G(u_1)=2$, we have
\begin{equation}\label{eq:3.15}
\left\{\begin{array}{ll}
  r+2t+(a+b-1)(4-t)=-ab-b^2+2b+5\ \ \ \ (t=2,3,4)\\[3pt]
  b=-2
   \end{array}
 \right.
\end{equation}
For $d_G(u_1)=a+b-1$, we have
\begin{equation}\label{eq:3.16}
\left\{\begin{array}{ll}
  r+2t+(a+b-1)(4-t)=-ab-b^2+2b+5\ \ \ \ (t=0,1,2)\\[3pt]
  d_G(u_2)=-ab-b^2-a+b+3
   \end{array}
 \right.
\end{equation}
Note that $d_G(u_2)\in\{2,a+b-1\}$. In view of (\ref{eq:3.15}),
since a and b are bother integers satisfying $a+b-1>5$, it's routine
to check that only $r=3$ with $t=3$ holds. It implies that
$(a,b)=(9,-2)$. Hence, $d_G(u_2)=a+b-1=6$. Applying (\ref{eq:1.1})
at $u_2$ yields $d_G(u_3)=10\notin\{2,a+b-1=6\}$, a contradiction.
In view (\ref{eq:3.16}), we get no integer solution satisfying
$a+b-1>5$, a contradiction. \vspace{2mm}

\textbf{Case 2}. $u$ is a non-attached vertex, i.e. $d_G(u)=5$. In
this case, we have $a+b-1\geq3$.\vspace{2mm}

If $d_G(u_1)=a+b-1$, then we apply (\ref{eq:1.1}) at $u_1$ to get
$d_G(u_2)=-ab-b^2+2b-3\in\{2,a+b-1\}$. Note that for $d_G(u_2)=2$,
we get $-b(a+b-2)=5$. As $a+b-1\geq3$ and $b$ is an integer, we have
$b=-1,a+b-2=5$, i.e., $(a,b)=(8,-1)$. If $d_G(u_1)=2$, then we apply
(\ref{eq:1.1}) at $u_1$ to get $d_G(u_2)=2a+b-9$. Note that
\[\label{eq:3.13}
d_G(u_2)\in\{2,a+b-1\}.
\]

Applying (\ref{eq:1.1}) at $u$, we have
\begin{equation}\label{eq:3.19}
\sum_{w\in N_G(u)}d_G(w)=5a+b-25.
\end{equation}

Together with (\ref{3-03}), we have $r+2t+(a+b-1)(4-t)=5a+b-25$. For
$d_G(u_1)=2$, we have
\begin{equation}\label{eq:3.20}
 \left\{\begin{array}{ll}
   r+2t+(a+b-1)(4-t)=5a+b-25\ \ \ \ \ (t=2,3,4)\\[3pt]
   d_G(u_2)=2a+b-9
   \end{array}
 \right.
\end{equation}
For $d_G(u_1)=a+b-1$, we have
\begin{equation}\label{eq:3.21}
\left\{\begin{array}{ll}
  r+2t+(a+b-1)(4-t)=5a+b-25\ \ \ \ \ (t=0,1,2)\\[3pt]
  d_G(u_2)=-ab-b^2+2b-3
   \end{array}
 \right.
\end{equation}

In view of (\ref{eq:3.13}) and (\ref{eq:3.20}), note that a and b
are both integers satisfying $a+b-1\geq3$. It's routine to check
that only $d_G(u_2)=a+b-1$ with $r\in\{3,a+b-1\}$ holds. If $r=3$,
we have $(a,b)=(8,-4)$. Then $d_G(u_2)=a+b-1=3$. Applying
(\ref{eq:1.1}) at $u_2$, we have $d_G(u_3)=8\notin\{2,5,a+b-1=3\}$,
a contradiction. If $r=a+b-1$, we have $a=8$ and $(b+5)(t-4)=0$,
which implies that $t=4$, thus $b\in\{-1,-2,-3,-4\}$ since
$a+b-1\geq3$ and $b\le -1$. For $b=-1$, we have $(a,b)=(8,-1)$.
Since $d_G(u_1)=2,d_G(u_2)=6$, applying (\ref{eq:1.1}) at $u_2$,
it's simple to verify that
$d_G(u_3)=5\notin\{d_{\widetilde{G}}(u_3)=2,a+b-1=6\}$, it's a
contradiction; for $b\in \{-2,-3,-4\}$, we will also get a
contradiction by the similar proof as above.
%If $b=-2$, then $(a,b)=(8,-2)$. As
%$d_G(u_1)=2,d_G(u_2)=a+b-1=5$, applying (\ref{eq:1.1}) at $u_2$,
%it's simple to verify that
%$d_G(u_3)=8\notin\{d_{\widetilde{G}}(u_3)=2,a+b-1=5\}$, a
%contradiction; If $b=-3$, then $(a,b)=(8,-3)$, since
%$d_G(u_1)=2,d_G(u_2)=a+b-1=4$. Applying (\ref{eq:1.1}) at $u_2$,
%it's simple to verify that
%$d_G(u_3)=9\notin\{d_{\widetilde{G}}(u_3)=2,a+b-1=4\}$, a
%contradiction; If $b=-4$, then $(a,b)=(8,-4)$, since
%$d_G(u_1)=2,d_G(u_2)=a+b-1=3$. Applying (\ref{eq:1.1}) at $u_2$,
%it's simple to verify that
%$d_G(u_3)=8\notin\{d_{\widetilde{G}}(u_3)=2,a+b-1=3\}$, a
%contradiction.

In view of (\ref{eq:3.13}) and (\ref{eq:3.21}), if $d_G(u_2)=2$,
then $(a,b)=(8,-1)$. For $r=2$, we have $t=3$; for $r=3$, we have
$t=\frac{13}{4}$; for $r=a+b-1$, we have $t=4$. Each is impossible
since $t\in\{0,1,2\}$. If $d_G(u_2)=a+b-1$, it's routine to check
that there is no integer solution such that $a+b-1\geq3$ since
$r\in\{2,3,a+b-1\}$ and $t=0,1,2$.
\end{proof}
\begin{figure}%[h!]
\begin{center}
  %\psfrag{a}{$G_{24}\in\mathscr{G}_{6,-1}$}
%\psfrag{b}{$G_{25}\in\mathscr{G}_{6,-2}$}\psfrag{c}{$G_{26}\in\mathscr{G}_{6,-2}$}
 \includegraphics[width=50mm]{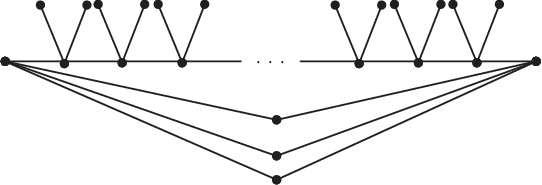}\\
  \caption{Graphs $G_{32}\in \mathscr{G}_{7,-2}.$ }
\end{center}
\end{figure}
\begin{prop}
Let $G\in \mathscr{G}_{a,b}$ be a tricyclic graph with pendants with
$\widetilde{G}=T_{12}$. Then $G\cong G_{32}\in\mathscr{G}_{7,-2},$
where $G_{29}$ is depicted in Fig. 4.
\end{prop}
\begin{proof}
Note that $\widetilde{G}=T_{12}$ (see Fig. 1), hence $\widetilde{G}$
consists of four internal paths connecting $u$ and $v$:
$P_{k_1+1}=us_1\ldots s_{k_1}(s_{k_1}=v),P_{k_2+1}=ut_1\ldots
t_{k_2}(t_{k_2}=v), P_{k_3+1}=ux_1\ldots
x_{k_3}(x_{k_3}=v),P_{k_4+1}=uy_1\ldots y_{k_4}(y_{k_4}=v)$. Without
loss of generality, we assume $k_4 \ge k_3\ge k_2 \geq k_1$. If
$k_1=1$, then $N_{\widetilde{G}}(u)=\{v, t_1,x_1,y_1\}$; if
$k_1\geq2$, then $N_{\widetilde{G}}(u)=\{s_1,t_1,x_1,y_1\}$ with
$s_1 \not= v$. Obviously, $k_4 \ge k_3\ge k_2\geq2$, then
$d_G(t_1),d_G(x_1),d_G(y_1)\in\{2,a+b-1\}$. Let
\[\label{3-02}
t=|\{x: d_G(x)=2,\ \ x\in \{t_1, x_1, y_1\}\}|,
\]
hence $0\le t\le 3.$ Denote
$N_{\widetilde{G}}(u)=\{w,t_1,x_1,y_1\}$, where  $w=v$ if $k_1=1$;
$w=s_1\not=v$ if $k_1\geq2$. Thus, $d_G(w)=r\in\{2,4,a+b-1\}$.
\vspace{2mm}

{\bf Case 1.}\ $d_G(u)=a+b-1>4$. In this case, since $b\leq-1$,  we have $a\geq7$. \vspace{2mm}

If $d_G(t_1)=2$, applying (\ref{eq:1.1}) at $t_1$, we have
$d_G(t_2)=a-3\in\{2,4,a+b-1\}$. Since $a\geq7$, we get that $a-3\neq
2$. For $d_G(t_2)=4$, we get $(a,b)=(7,-1)$; for $d_G(t_2)=a+b-1$,
we get $b=-2$. If $d_G(t_1)=a+b-1$, applying (\ref{eq:1.1}) at
$t_1$, we have $d_G(t_2)=-ab-b^2-a+b+3$.

Applying (\ref{eq:1.1}) at $u$ yields
\begin{equation}\label{eq:3.30}
d_G(w)+d_G(t_1)+d_G(x_1)+d_G(y_1)=-ab-b^2+2b+4.
\end{equation}

By (\ref{3-02}) we have
\begin{equation}\label{eq:3.31}
r+2t+(a+b-1)(3-t)=-ab-b^2+2b+4, \ \ \ \ \ (t=0, 1, 2, 3)
\end{equation}

$\bullet$ $t=0$. (\ref{eq:3.31}) gives $ab+b^2+3a+b=7-r$. Since
$d_G(t_1)=a+b-1$, we have
\begin{equation*}\label{eq:3.32}
\text{$ab+b^2+3a+b=7-r$  \ \  and \ \
    $d_G(t_2)= -ab-b^2-a+b+3$}
\end{equation*}
with $d_G(t_2)\in\{2,4,a+b-1\}$. It's routine to check that there is
no integer solution such that $a+b-1>4$, a contradiction.
%If $d_G(t_2)=2$, then
%(\ref{eq:3.32}) gives $2(a+b)=6-r$, a contradiction to the fact that
%$a+b-1>4$. If $d_G(t_2)=4$, then (\ref{eq:3.32}) gives $2a+2b=8-r$.
%It's impossible since $a+b-1>4$. If $d_G(t_2)=a+b-1$,
%(\ref{eq:3.32}) gives $a+b=3+r$, so $r\neq a+b-1$; on the other
%hand, (\ref{eq:3.32}) implies $rb+2+4r+b=0$, i.e.
%$b=\frac{-2-4r}{r+1}$. Since $b$ is an integer, it's simple to
%verify that $r\notin \{2,4\}$, a contradiction.

$\bullet$ $t=1$. (\ref{eq:3.31}) gives $ab+b^2+2a=4-r$. Without loss
of generality, we assume $d_G(t_1)=2$, then $(a,b)=(7,-1)$ or
$b=-2$. It's easy to check that both are impossible since
$r\in\{2,4,a+b-1\}$.

$\bullet$ $t=2$. (\ref{eq:3.31}) gives $ab+b^2+a-b=1-r$. Without
loss of generality, we assume $d_G(t_1)=2$, then $(a,b)=(7,-1)$ or
$b=-2$. If $(a,b)=(7,-1)$, then $r=-1$, a contradiction; if $b=-2$
with $d_G(t_2)=a+b-1$, then $a=5+r$. It's easy to see that only
$r=4$ holds. So we have $v\in N_G(u)$ and $d_G(v)=4$, moreover,
$(a,b)=(9,-2)$ and $d_G(t_2)=6$. Applying (\ref{eq:1.1}) at $t_2$,
we have $d_G(t_3)=10\notin\{2,4,a+b-1=6\}$, a contradiction.

$\bullet$ $t=3$. (\ref{eq:3.31}) gives $ab+b^2-2b+2+r=0$. Note that
$d_G(t_1)=2$, hence $(a,b)=(7,-1)$ or $b=-2$.

If $(a,b)=(7,-1)$, then $r=2$. Therefore, we have
$d_G(s_1)=d_G(t_1)=d_G(x_1)=d_G(y_1)=2$ together with $t=3$.
Applying (\ref{eq:1.1}) at $s_1$ yields $d_G(s_2)=4$, which implies
$s_2=v$ with $d_G(v)=4$. Similarly, $t_2=x_2=y_2=v$. Thus,
$N_G(v)=\{s_1,t_1,x_1,y_1\}$. We may check that this is impossible
by applying (\ref{eq:1.1}) at $v$.

If $b=-2$ with $d_G(t_2)=a+b-1$, then $2a=10+r$. Since $a+b-1>4$ and
$a$ is an integer, it's easy to verify that $r\notin \{2,4,a+b-1\}$,
a contradiction. \vspace{2mm}

{\bf Case 2.}\ $d_G(u)=4$. In this case, $a+b-1\geq 3$. \vspace{2mm}

If $d_G(t_1)=2$, applying (\ref{eq:1.1}) at $t_1$, we have
$d_G(t_2)=2a+b-8$. If $d_G(t_1)=a+b-1$, applying (\ref{eq:1.1}) at
$t_1$, we have $d_G(t_2)=-ab-b^2+2b-2$.

Applying (\ref{eq:1.1}) at $u$ yields
\begin{equation}\label{eq:3.33}
d_G(w)+d_G(t_1)+d_G(x_1)+d_G(y_1)=4a+b-16
\end{equation}

From (\ref{3-02}) we have
\begin{equation}\label{eq:3.34}
r+2t+(a+b-1)(3-t)=4a+b-16, \ \ \ \ (t=0, 1, 2, 3.)
\end{equation}

$\bullet$ $t=0$. (\ref{eq:3.34}) gives $a-2b=13+r$. Since
$d_G(t_1)=a+b-1$, we have
\begin{equation*}\label{eq:3.35}
\text{$a-2b=13+r$  \ \  and \ \
    $d_G(t_2)=-ab-b^2+2b-2$}.
\end{equation*}

Note that $d_G(t_2)\in\{2,4,a+b-1\}$. It's routine to check that
there is no integer solution such that $a+b-1\geq3$ since
$r\in\{2,4,a+b-1\}$.

$\bullet$ $t=1$. (\ref{eq:3.34}) gives $2a-b=16+r$. Without
loss of generality, we assume $d_G(t_1)=a+b-1$, then
\begin{equation*}\label{eq:3.36}
\text{$2a-b=16+r$  \ \  and \ \
    $d_G(t_2)=-ab-b^2+2b-2$}.
\end{equation*}

Note that $d_G(t_2)\in\{2,4,a+b-1\}$. It's routine to check that
there is no integer solution such that $a+b-1\geq3$ since
$r\in\{2,4,a+b-1\}$.

$\bullet$ $t=2$. (\ref{eq:3.34}) gives that $3a=19+r$. Since $a$ is
an integer, $r\neq4$. Hence, $r\in\{2,a+b-1\}$. Without loss of
generality, we assume $d_G(t_1)=d_G(x_1)=2,d_G(y_1)=a+b-1$, then
\begin{equation*}\label{eq:3.36}
\text{$3a=19+r$  \ \  and \ \
    $d_G(t_2)=2a+b-8$}.
\end{equation*}

Note that $d_G(t_2)\in\{2,4,a+b-1\}$. It's routine to check that
only $r=2$ holds, which implies that $a=7$, moreover, $d_G(s_1)=2$.
For $d_G(t_2)=2$, we have $b=-4$. However, $a+b-1=2<3$, a
contradiction. For $d_G(t_2)=4$, we have $b=-2$. For
$d_G(t_2)=a+b-1$, we have $b\in\{-1,-2,-3\}$.

First we consider $(a,b)=(7,-2)$. In this subcase, $a+b-1=4$ and
$d_G(s_2)=d_G(x_2)=d_G(t_2)=4$. If $\{s_2,t_2,x_2\}$ contains a
member, say $s_2$, such that $s_2\neq v$, then applying
(\ref{eq:1.1}) at $s_2$ yields $d_G(s_3)=6\notin\{2,4\}$, a
contradiction. Therefore, we get that $s_2=x_2=t_2=v$. Note that
$d_G(y_1)=a+b-1=4$. It's easy to see that $y_1\not=v$ since G has
pendant vertices. Applying (\ref{eq:1.1}) at $y_1$, we have
$d_G(y_2)=4$. Continue the process, we may finally obtain that
$d_G(y_i)=4$ for all $1\leq i\leq k_4$, where $k_4\geq2$. Thus, we
get the graph $G\cong G_{32}\in\mathscr{G}_{7,-2}$; see Fig. 4.

Next, we consider $b\in\{-1,-3\}$.

If $b=-1$, then $(a,b)=(7,-1).$ In this subcase, $d_G(y_1)=a+b-1=5$,
and $d_G(s_2)=d_G(t_2)=d_G(x_2)=5$. If $\{s_2,t_2,x_2\}$ contains a
member, say $s_2$, such that $s_2\neq v$, then applying
(\ref{eq:1.1}) at $s_2$ yields $d_G(s_3)=4$. It implies that $s_3=v$
with $d_G(v)=4$. Therefore, $t_2\not=v,x_2\not=v$ and $t_3=x_3=v$.
It's easy to see that this is impossible by applying (\ref{eq:1.1})
at $v$. Hence, we have $s_2=t_2=x_2=v$ with $d_G(v)=5$. Note that
$d_G(y_1)=5$. We apply (\ref{eq:1.1}) at $y_1$ to get $d_G(y_2)=2$,
moreover, $d_G(y_3)=4\notin\{2,5\}$, a contradiction.

If $b=-3$, then $(a,b)=(7,-3)$. In this subcase, $d_G(y_1)=a+b-1=3$,
and $d_G(s_2)=d_G(t_2)=d_G(x_2)=3$. Applying (\ref{eq:1.1}) at $s_2$
yields $d_G(s_3)=6\notin\{2,4,a+b-1=3\}$, a contradiction.

$\bullet$ $t=3$. (\ref{eq:3.34}) gives that $4a+b=22+r$. As
$d_G(t_1)=2$, we have
\begin{equation}\label{eq:3.37}
\text{$4a+b=22+r$  \ \  and \ \
    $d_G(t_2)=2a+b-8$}.
\end{equation}

Note that $d_G(t_2)\in\{2,4,a+b-1\}$. If $d_G(t_2)=2$, then
(\ref{eq:3.37}) gives that $2a=12+r$. It's routine to check that
there is no integer solution such that $a+b-1\geq3$, a
contradiction.
%If $r\in\{2,a+b-1\}$, so $(a,b)=(7,-4)$, a contradiction since $a+b-1=2<3$. If $r=4$, then $(a,b)=(8,-6)$, a contradiction since $a+b-1=1<3$.

If $d_G(t_2)=4$, then (\ref{eq:3.37}) gives that $2a=10+r$. For
$r=2$, we have $b=0$, a contradiction to $b\leq -1$. For
$r\in\{4,a+b-1\}$, we have $(a,b)=(7,-2)$.

If $d_G(t_2)=a+b-1$, then (\ref{eq:3.37}) gives that $a=7$, this
implies that $a+b-1=r$. Note that $a+b-1\geq3$, then $r\not=2$. For
$r=4$, we have $(a,b)=(7,-2)$. For $r=a+b-1$, we have
$b\in\{-1,-2,-3\}$ since $b\leq-1$. By a similar proof as in the
discussion of $t=2$, we may also get $G\cong
G_{32}\in\mathscr{G}_{7,-2}$.

This completes the proof.
\end{proof}

By a similar discussion as in the proof of Lemma 3.5, we can show the next lemma. We omit its procedure.
\begin{prop}
Let $G$ be a tricyclic graph with pendants. If $\widetilde{G}=T_{13}$ (see Fig. 1), then $G\notin\mathscr{G}_{a,b}$.
\end{prop}
\begin{figure}[h!]
\begin{center}
 \psfrag{d}{$G_{36}\in\mathscr{G}_{7,-1}$} \psfrag{a}{$G_{33}\in\mathscr{G}_{6,-1}$}
\psfrag{b}{$G_{34}\in\mathscr{G}_{6,-2}$}\psfrag{c}{$G_{35}\in\mathscr{G}_{6,-2}$}
 \includegraphics[width=150mm]{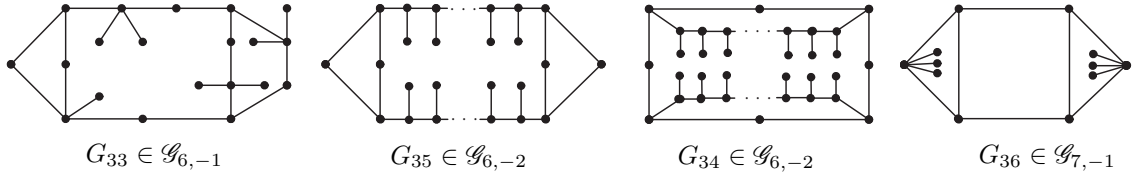}\\
  \caption{Graphs $G_{33}, \, G_{34}, \, G_{35},\, G_{36}.$ }
\end{center}
\end{figure}
\begin{prop}
Let $G\in \mathscr{G}_{a,b}$ be a tricyclic graph with pendants
satisfying $\widetilde{G}=T_{14}$. Then $G\cong
G_{33}\in\mathscr{G}_{6,-1}$ or, $G\cong
G_{34}\in\mathscr{G}_{6,-2}$ or, $G\cong
G_{35}\in\mathscr{G}_{6,-2}$ or, $G\cong
G_{36}\in\mathscr{G}_{7,-1}$, where $G_{33}, \, G_{34},\, G_{35},\,
G_{36}$ are depicted in Fig. 5.
\end{prop}
\begin{proof}
Note that $\widetilde{G}=T_{14}$; see Fig. 1, hence $\widetilde{G}$
contains two internal paths $P_{k_1+1}=u_1s_1s_2\ldots
s_{k_1}(s_{k_1}=u_2), P_{k_2+1}=u_1t_1t_2\ldots
t_{k_2}(t_{k_2}=u_2)$ connecting $u_1,\,u_2$, one internal path
$P_{k_3+1}=u_1x_1x_2\ldots x_{k_3}(x_{k_3}=v_1)$ connecting
$u_1,\,v_1$, one internal path $P_{k_4+1}=u_2y_1y_2\ldots
y_{k_4}(y_{k_4}=v_2)$ connecting $u_2,v_2$, two internal paths
$P_{k_5+1}=v_1z_1z_2\ldots z_{k_5}(z_{k_5}=v_2),
P_{k_6+1}=v_1q_1q_2\ldots q_{k_6}(q_{k_6}=v_2)$ connecting
$v_1,v_2$. Without loss of generality, We assume $k_2\ge 2$. Note
that $N_{\widetilde{G}}(u_1)=\{s_1,t_1,x_1\}$. If $k_3=1$, then
$x_1=v_1$; if $k_3\geq2$, then $x_1\not=v_1$. Set $m=|\{x: d_G(x)=3,
x\in \{s_1,t_1,x_1\}\}|.$ According to the structure of
$\widetilde{G}$, we have $m=0,1,2.$ We consider the following two
possible cases.

\vspace{2mm} \textbf{Case 1}.  There exists at least one vertex in
$\{u_1, u_2, v_1, v_2\}$, say $u_1$, such that it is an attached
vertex. In this case, we have $a+b-1\geq3$.\vspace{2mm}

If $d_G(t_1)=2$, applying (\ref{eq:1.1}) at $t_1$ yields
$d_G(t_2)=2a+b-7$; if $d_G(t_1)=a+b-1$, applying (\ref{eq:1.1}) at
$t_1$ yields $d_G(t_2)=-ab-b^2+2b-1$.

Applying (\ref{eq:1.1}) at $u_1$, we get
\begin{equation}\label{eq:3.01}
d_G(s_1)+d_G(t_1)+d_G(x_1)=3a+b-9
\end{equation}

\vspace{2mm}\textbf{Subcase 1.1}. $m=0,$ or $1.$ In this case, we
have $\{d_G(s_1),d_G(t_1),d_G(x_1)\}=\{r,t\cdot 2,(2-t)\cdot
(a+b-1)\}$, where $t=0, 1, 2$ and $r\in\{2, 3, a+b-1\}$. It is easy
to see that $m=1$ if $r=3$ and 0 otherwise.\vspace{2mm}

In view of (\ref{eq:3.01}), we have
\begin{equation}\label{eq:3.02}
r+2t+(a+b-1)(2-t)=3a+b-9,\, \ \ \ \ t=0,1,2
\end{equation}

$\bullet$ $t=0$. Then (\ref{eq:3.02}) gives $a-b=7+r$. In this
subcase we may assume $d_G(t_1)=a+b-1$. Hence, we have
\begin{equation}\label{eq:3.41}
\text{$a-b=7+r$ \ \ and \ \
    $d_G(t_2)=-ab-b^2+2b-1$}.
  \end{equation}
Note that $d_G(t_2)\in\{2,3,a+b-1\},\, r\in\{2,3,a+b-1\}$.
%If $r\in\{2,a+b-1\}$, then (\ref{eq:3.41}) gives that $(a,b)=(6,-3)$, a contradiction to the fact $a+b-1\geq3$. If $r=3$,
It is routine to check that (\ref{eq:3.41}) has no integer solution
satisfying $a+b-1\geq3$, a contradiction.

$\bullet$ $t=1$. Then at least one of $\{s_1,t_1,x_1\}$ is of degree
two. We  first assume that at least one vertex in $\{s_1, t_1\}$ is
of degree 2. For convenience, let $d_G(t_1)=2$. On the other hand,
since $t=1$, (\ref{eq:3.02}) gives $2a=10+r$. Hence,
\begin{equation}\label{eq:3.03}
\text{$2a=10+r$  \ \  and \ \
    $d_G(t_2)=2a+b-7$}
\end{equation}
with $d_G(t_2)\in\{2,3,a+b-1\}$. It's routine to check that only
$r=2$ holds, which implies that $a=6$. If $d_G(t_2)=2$, then
(\ref{eq:3.03}) gives $b=-3,a+b-1=2<3$, a contradiction. If
$d_G(t_2)=3$, then (\ref{eq:3.03}) gives $b=-2$. Hence,
$(a,b)=(6,-2)$. If $d_G(t_2)=a+b-1$, then $b\in\{-1,-2\}$ since
$a+b-1\geq3$ and $b\leq-1$.

First we consider $(a,b)=(6,-1)$ with $d_G(t_2)=a+b-1=4$. Note that
$d_G(s_1)\in\{2,4\}$. If $d_G(s_1)=2$, then $d_G(x_1)=4$. Applying
(\ref{eq:1.1}) at $s_1$ yields $d_G(s_2)=4$. Next, we will show that
$s_2=t_2=u_2$. On the contrary, we suppose $s_2\not=t_2$. Applying
(\ref{eq:1.1}) at $s_2$, we have $d_G(s_3)=3$, which implies that
$s_3=u_2$ with $d_G(u_2)=3$. Hence, $t_2\not=u_2$ and $t_3=u_2$.
Then applying (\ref{eq:1.1}) at  $u_2$ yields $d_G(y_1)=0$, a
contradiction. Therefore, we get that $s_2=u_2$, similarly,
$t_2=u_2$. Thus, $d_G(u_2)=4$. Applying (\ref{eq:1.1}) at $u_2$
yields $d_G(y_1)=2$, moreover, $d_G(y_2)=4, d_G(y_3)=3$. Hence,
$y_3=v_2$. Note that $d_G(x_1)=4$. It's easy to see that
$x_1\not=v_1$. In fact, if $x_1=v_1$, then applying (\ref{eq:1.1})
at $v_1$ yields $d_G(z_1)+d_G(q_1)=3$. This is impossible.
Therefore, we get that $x_1\not=v_1$. We apply (\ref{eq:1.1}) at
$x_1$ to get $d_G(x_2)=2$, moreover, $d_G(x_3)=3$. So we have
$x_3=v_1$ and $d_G(v_1)=3$. Applying (\ref{eq:1.1}) at $v_1$, we
have $d_G(z_1)+d_G(q_1)=6$. Without loss of generality, we assume
$d_G(z_1)=2,d_G(q_1)=4$. We apply (\ref{eq:1.1}) at $z_1$ and $q_1$
respectively to get $d_G(z_2)=4,d_G(q_2)=2$. Moreover,
$d_G(z_3)=d_G(q_3)=3$. Hence, $z_3=q_3=v_2$. Thus, it's easy to
check that $G\cong G_{33}\in\mathscr{G}_{6,-1}$; See Fig. 5.

If $d_G(s_1)=4$, then $d_G(x_1)=2$. By a similar proof above, we may
also obtain that  $G\cong G_{33}\in\mathscr{G}_{6,-1}$.

Now we consider $(a,b)=(6,-2)$ with $d_G(t_2)=3=a+b-1$. If $t_2\neq
u_2$, then applying (\ref{eq:1.1}) at $t_2$ yields $d_G(t_3)=4
\notin\{2,3\}$, a contradiction. Hence, $t_2=u_2$. Note that
$d_G(s_1)\in\{2,3\}$. For $d_G(s_1)=2$, we may get $s_2=u_2$. By
(\ref{eq:1.1}), it is routine to check that
$d_G(x_1)=\cdots=d_G(x_{k_3})=3,d_G(y_1)=\cdots=d_G(y_{k_4})=3(
k_3\ge 1, k_4\ge 1$ and $k_3k_4\ge 1)$. Furthermore, $d_G(z_1)=2,
z_2=v_2; d_G(q_1)=2, q_2=v_2$. Therefore, we obtain that $G\cong
G_{34}\in \mathscr{G}_{6,-2};$ see Fig. 5.

For $d_G(s_1)=3$, we have $d_G(x_1)=2$ and $x_2=v_1$. By
(\ref{eq:1.1}), it is routine to check that
$d_G(s_1)=\cdots=d_G({s_{k_1}})=3(k_1\geq1)$. Similarly,
$d_G(y_1)=2$ and $y_2=v_2$; $d_G(z_1)=2,z_2=v_2;
d_G(q_1)=\cdots=d_G(q_{k_6})=3(k_6\geq1)$. Note that G has pendant
vertices, then $k_1k_6\not=1)$. Hence, we obtain that $G\cong
G_{35}\in \mathscr{G}_{6,-2};$ see Fig. 5.

Now we assume $d_G(t_1), d_G(s_1)\not=2.$ Hence, $d_G(x_1)=2$. By a
similar discussion as the former subcase, we have $r=2$. Together
with $t=1$, we have at least two of $\{s_1,t_1,x_1\}$ of degree 2, a
contradiction to the assumption.

$\bullet$ $t=2$. Then (\ref{eq:3.02}) gives $3a+b=13+r$. Without
loss of generality, we assume that $d_G(t_1)=2$. Hence,
\begin{equation}\label{eq:3.43}
\text{$3a+b=13+r$  \ \  and \ \
    $d_G(t_2)=2a+b-7$}.
\end{equation}
Note that $d_G(t_2)\in\{2,3,a+b-1\}$, $r\in\{2,3,a+b-1\}$.
%First consider $r=3$. If $d_G(t_2)=2$, then by (\ref{eq:3.43}) we have $(a,b)=(7,-5)$, a
%contradiction. Hence, $d_G(t_2)\in\{3,a+b-1\}$.
It's easy to check that $r\not=2,d_G(t_2)\not=2$. Furthermore, for
$d_G(t_2)=3$ with $r\in\{3,a+b-1\}$ or $d_G(t_2)=a+b-1$ with $r=3$,
we have $(a,b)=(6,-2)$; for $d_G(t_2)=a+b-1$ with $r=a+b-1$, we have
$a=6,b\in\{-1,\,-2\}$.

First we consider $(a,b)=(6,-1)$.  In this subcase, $r=a+b-1=4$ and
$\{d_G(s_1),d_G(t_1),d_G(x_1)\}=\{2,2,4\}$. By a similar discussion
as in the proof of the case $t=1$ and $(a,b)=(6,-1)$, we can also
obtain the graph $G\cong G_{33}\in \mathscr{G}_{6,-1}$.

Now  we consider $(a,b)=(6,-2)$. In this subcase, $r=a+b-1=3$ and
$\{d_G(s_1),d_G(t_1),d_G(x_1)\}=\{2,2,3\}$. By a similar discussion
as in the proof of the case $t=1$ and $(a,b)=(6,-2)$, we can also
obtain the graph $G\cong G_{34},G_{35}\in \mathscr{G}_{6,-2}$. \vspace{2mm}

\vspace{2mm}\textbf{Subcase 1.2}.  $m=2$, i.e.,
$\{d_G(s_1),d_G(t_1),d_G(x_1)\}=\{2,3,3\}$ or
$\{d_G(s_1),d_G(t_1),d_G(x_1)\}=\{a+b-1,3,3\}$.\vspace{2mm}

First we consider $\{d_G(s_1),d_G(t_1),d_G(x_1)\}=\{2,3,3\}$. In
this subcase, we let $d_G(s_1)=2,d_G(t_1)=d_G(x_1)=3$. Thus,
$t_1=u_2,x_1=v_1$ with $d_G(u_2)=d_G(v_1)=3$.

In view of (\ref{eq:3.01}), we have $2+3+3=3a+b-9$, i.e., $3a+b=17$.
Applying (\ref{eq:1.1}) at $s_1$ yields
$d_G(s_2)=2a+b-7\in\{2,3,a+b-1\}$. It's routine to check that only
$d_G(s_2)=a+b-1$ holds, which implies $(a,b)=(6,-1)$. Note that
$s_2\not=u_2$. Applying (\ref{eq:1.1}) at $s_2$ yields $d_G(s_3)=3$,
thus, $s_3=u_2$. Applying (\ref{eq:1.1}) at $u_2$, we get that
$d_G(y_1)=1$, a contradiction.

Now consider that $\{d_G(s_1),d_G(t_1),d_G(x_1)\}=\{a+b-1,3,3\}$. In
this subcase, we let $d_G(s_1)=a+b-1,d_(t_1)=d_G(x_1)=3$. Thus,
$t_1=u_2,x_1=v_1$ with $d_G(u_2)=d_G(v_1)=3$.

In view of (\ref{eq:3.01}), we have $a+b-1+3+3=3a+b-9$, which gives
that $a=7$. Applying (\ref{eq:1.1}) at $s_1$ yields
$d_G(s_2)=-ab-b^2+2b-1\in\{2,3,a+b-1\}$. It's routine to check that
only $d_G(s_2)=3$ holds, which implies $(a,b)=(7,-1)$ and $s_2=u_2$.
Applying (\ref{eq:1.1}) at $u_2$, we get $v_2\in N_G(u_2)$ and
$d_G(v_2)=3$. Applying (\ref{eq:1.1}) at $v_1$ yields
$d_G(z_1)+d_G(q_1)=8$. Without loss of generality, we assume
$d_G(z_1)=3,d_G(q_1)=5$. Thus, $z_1=v_2$. Applying (\ref{eq:1.1}) at
$q_1$, we have $d_G(q_2)=3$. Hence, $q_2=v_2$. It's simple to check
that $G\cong G_{36}\in\mathscr{G}_{7,-1}$; See Fig. 5.

\vspace{2mm} \textbf{Case 2}. Each of the vertices in
$\{u_1,\,u_2,\,v_1,\,v_2\}$ is an attached-vertex. That is,
$d_G(u_1)=d_G(u_2)=d_G(v_1)=d_G(v_2)=a+b-1>3$. In this case, we have
$\{d_G(s_1), d_G(t_1), d_G(x_1)\}=\{t\cdot 2, (3-t)\cdot (a+b-1)\},$
where $0\leq t\leq3$.\vspace{2mm}

If $d_G(t_1)=2$, applying (\ref{eq:1.1}) at $t_1$, we have
$d_G(t_2)=a-3\in\{2,a+b-1\}$. Since $a\geq6$, we obtain $a-3\neq 2$.
For $d_G(t_2)=a+b-1$, we have $b=-2$. If $d_G(t_1)=a+b-1$, applying
(\ref{eq:1.1}) at $t_1$, we have $d_G(t_2)=-ab-b^2-a+b+3$.

Applying Lemma \ref{lem4}(ii) at $u_1$, we get
$d_G(s_1)+d_G(t_1)+d_G(x_1)=-ab-b^2+2b+3$, which gives
\begin{equation}\label{eq:3.38}
2t+(a+b-1)(3-t)=-ab-b^2+2b+3,\, \ \ \ t=0,1,2,3.
\end{equation}

$\bullet$ $t=0$. (\ref{eq:3.38}) gives $ab+b^2+3a+b=6$. In this
subcase, we have $d_G(t_1)=a+b-1$. Hence,
\begin{equation}\notag
\text{$ab+b^2+3a+b=6$  \ \ and \ \
     $d_G(t_2)=-ab-b^2-a+b+3$}
\end{equation}
with $d_G(t_2)\in\{2,a+b-1\}$. It's routine to check that there is
no integer solution such that $a+b-1>3$.

$\bullet$ $t=1$. (\ref{eq:3.38}) gives $ab+b^2+2a=3$. If
$d_G(t_1)=2$, then $b=-2$. However, it does not satisfy the
equation. Therefore, $d_G(t_1)=a+b-1$. Hence,
\begin{equation}\label{eq:3.39}
\text{$ab+b^2+2a=3$  \ \ and \ \
     $d_G(t_2)=-ab-b^2-a+b+3$}.
\end{equation}
By (\ref{eq:3.39}), we get that $d_G(t_2)=a+b\in\{2,a+b-1\}$. It's
impossible since $a+b-1>3$.

$\bullet$ $t=2,3$. Without loss of generality we assume
$d_G(t_1)=2$. Then we have $b=-2$. For $t=2$, (\ref{eq:3.38}) gives
$ab+b^2+a-b=0$; for $t=3$, (\ref{eq:3.38}) gives $ab+b^2-2b+3=0$.
It's easy to check that both are impossible since $b=-2$.

Thus, we complete the proof.
\end{proof}

\begin{prop}
Let $G\in \mathscr{G}_{a,b}$ be a tricyclic graph with
$\widetilde{G}=T_{15}$. Then $G\cong G_{34}\in\mathscr{G}_{8,-2},$
or $G\cong G_{35}\in\mathscr{G}_{7,-2},$ or
$G_{36}\in\mathscr{G}_{7,-1},$ or $G\cong
G_{37}\in\mathscr{G}_{a,b}$ with $\frac{6-3a}{b}=k\geq6,$ or $G\cong
G_{38}\in\mathscr{G}_{8,-2},$  or $G\cong
G_{39}\in\mathscr{G}_{6,-1},$  or $G\cong
G_{40}\in\mathscr{G}_{6,-1},$ where $G_{34}, G_{35}, G_{36}, G_{37}$
and $G_{38}, G_{39}, G_{40}$ are depicted in Fig. 5.
\end{prop}
\begin{figure}[h!]
\begin{center}
\psfrag{a}{$k$}\psfrag{x}{$G_{34}\in\mathscr{G}_{8,-2}$}
\psfrag{y}{$G_{38}\in\mathscr{G}_{7,-2}$}\psfrag{z}{$G_{37}\in\mathscr{G}_{8,-2}$}
\psfrag{w}{$G_{39}\in\mathscr{G}_{7,-1}$}\psfrag{t}{$G_{40}\in\mathscr{G}_{6,-1}$}
\psfrag{n}{$G_{41}\in\mathscr{G}_{6,-1}$}\psfrag{m}{$G_{42}\in\mathscr{G}_{a,b}$}
\includegraphics[width=140mm]{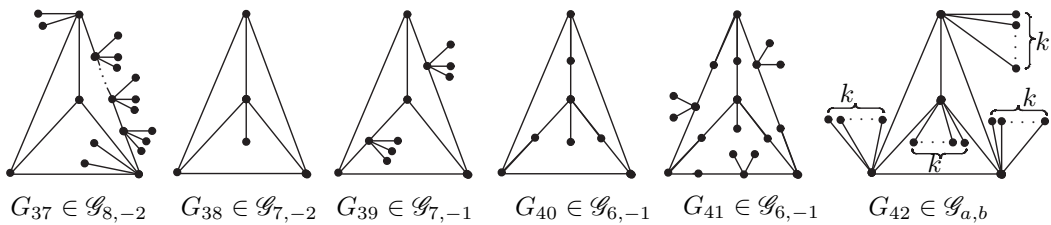}\\
\caption{Graphs $G_{37}, G_{38}, G_{39}, G_{40}, G_{41}$ and $G_{42}$.}
\end{center}
\end{figure}
\begin{proof}
Note that $\widetilde{G}=T_{15}$ (see Fig. 1), then
$d_{\widetilde{G}}(v_1)=d_{\widetilde{G}}(v_2)=d_{\widetilde{G}}(v_3)=d_{\widetilde{G}}(v_4)=3$
and $\widetilde{G}$ consists of six internal paths:
$P_{k_1}=v_2u_1\ldots u_{k_1} (u_{k_1}=v_3), P_{k_2}=v_2w_1\ldots
w_{k_2} (w_{k_2}=v_4), P_{k_3}=v_3s_1\ldots s_{k_3} (s_{k_3}=v_4),
P_{k_4}=v_1x_1\ldots x_{k_4} (x_{k_4}=v_2), P_{k_5}=v_1y_1\ldots
y_{k_5} (y_{k_5}=v_3),P_{k_6}=v_1z_1\ldots z_{k_6} (z_{k_6}=v_4)$.

\vspace{2mm} {\bf Case 1.}\ There exists at least one vertex in
$\{v_1,\,v_2,\,v_3,\,v_4\}$, say $v_1$, such that it is a
non-attached vertex. That is, $d_G(v_1)=3$. In this case,
$a+b-1\geq3$. Let $N_{\widetilde{G}}(v_1)=\{r_1, r_2, r_3\}$, then
$\{d_G(r_1), d_G(r_2), d_G(r_3)\}=\{t_1\cdot 2, t_2\cdot 3, t_3\cdot
(a+b-1)\},$ where $t_1+t_2+t_3=3$ with $0\leq t_i\leq3(i=1,2,3)$.

\vspace{2mm} \textbf{Subcase 1.1}. $t_1=0$, that is $\{d_G(r_1),
d_G(r_2), d_G(r_3)\}=\{t\cdot 3,\, (3-t)\cdot (a+b-1)\}$, where
$0\leq t\leq3$. \vspace{2mm}

Applying (\ref{eq:1.1}) at $v_1$ yields
$d_G(r_1)+d_G(r_2)+d_G(r_3)=3a+b-9$. Hence we have
\begin{equation}\label{eq:3.44}
   3t+(a+b-1)(3-t)=3a+b-9,\ \ \ \ (t=0, 1, 2, 3).
\end{equation}

\medskip\noindent
$\bullet$ $t=0$. Then (\ref{eq:3.44}) gives that $b=-3$. First we
assume that $v_2,v_3,v_4\in N_G(v_1)$, then
$d_G(v_2)=d_G(v_3)=d_G(v_4)=a+b-1>3$. Applying (\ref{eq:1.1}) at
$v_2$ yields that $d_G(u_1)+d_G(w_1)=-ab-b^2+2b$, whence we may
check that $d_G(u_1),d_G(w_1)$ will not be in $\{2,a+b-1\}$, a
contradiction.

So we assume that at least one member in $\{v_2,v_3,v_4\}$, say
$v_4$, is not in $N_G(v_1)$. Thus, $d_G(z_1)=a+b-1>2$. Applying
Lemma (1.8)(ii) at $z_1$ yields
$d_G(z_2)=-ab-b^2+2b-1\in\{2,3,a+b-1\}$. It's routine to check that
there is no integer solution such that $a+b-1>2$, a contradiction.
%If $d_G(z_2)\in\{2,a+b-1\}$, then $(a,b)=(6,-3)$, a contradiction. If $d_G(z_2)=3$, then $a=\frac{19}{3}$, a contradiction.

\medskip\noindent
$\bullet$ $t=1$. Without loss of generality, we suppose $d_G(v_2)=3$
with $v_2\in N_G(v_1)$. By (\ref{eq:3.44}), we have $a-b=10$.

If $v_3,v_4\in N_G(v_1)$, thus, $d_G(v_3)=d_G(v_4)=a+b-1>3$.
Applying (\ref{eq:1.1}) at $v_2$ yields that
$d_G(u_1)+d_G(w_1)=3a+b-12$. Notice that $d_G(u_1),d_G(w_1)\in\{2,
a+b-1\}$, it is easy to check that only $d_G(u_1)=d_G(w_1)=a+b-1$ is
true, which also implies that $a-b=10$.

We first assume $u_1=v_3$ and $w_1=v_4$. Then applying
(\ref{eq:1.1}) at $v_3$ yields
$d_G(s_1)=-ab-b^2+2b-3\in\{2,a+b-1\}$. Note that $a-b=10$, for
$d_G(s_1)=2$, there is no integer solution; for $d_G(s_1)=a+b-1$, we
get that $(a,b)=(8,-2)$ and $a+b-1=5$. It's easy to check that
$d_G(s_i)=5$ for $i=1,2,\ldots,k_3$, where $k_3\geq1$. Thus, we get
$G\cong G_{37}\in\mathscr{G}_{8,-2}$; see Fig. 6.

Now we assume, without loss of generality, that $w_1\neq v_4$. Then
applying (\ref{eq:1.1}) at $w_1$ yields $d_G(w_2)=-ab-b^2+2b-1\in
\{2,a+b-1\}$. It is easy to check that this is impossible.

If $\{v_3,v_4\}$ contains a member, say $v_4$, not in $N_G(v_1)$.
Thus, $d_G(z_1)=a+b-1$ with $z_1\not=v_4$. Applying (\ref{eq:1.1})
at $z_1$ yields $d_G(z_2)=-ab-b^2+2b-1\in\{2,3,a+b-1\}$. It is
routine to check that this is impossible.

\medskip\noindent
$\bullet$ $t=2$.  Without loss of generality, we assume
$d_G(v_2)=d_G(v_3)=3$ with $v_2,v_3\in N_G(v_1)$. By (\ref{eq:3.44})
we get $a=7$. First, we consider $v_4\in N_G(v_1)$, then
$d_G(v_4)=a+b-1>3$. Applying (\ref{eq:1.1}) at $v_2$ gives
$d_G(u_1)+d_G(w_1)=3a+b-12$. Note that $d_G(u_1)\in \{2,3,a+b-1\},
d_G(w_1)\in\{2,a+b-1\}$. We may check that only $d_G(u_1)=3$ with
$u_1=v_3$ and $d_G(w_1)=a+b-1$ is true, which also implies that
$a=7$.

If $w_1\neq v_4$, applying (\ref{eq:1.1}) at $w_1$, we have that
$d_G(w_2)=-ab-b^2+2b-1\notin\{2,a+b-1\}$, a contradiction. So we
have $w_1=v_4$. Similarly, by applying (\ref{eq:1.1}) at $v_3$, we
may get $d_G(s_1)=a+b-1$ and $s_1=v_4$. Thus, we get the graph
$G\cong G_{39}$. We apply Lemma 1.8(ii) at $v_4$ to get
$3+3+3=-ab-b^2+2b+3$, together with $a=7$, we have $b=-2$. Hence
$G\cong G_{38}\in\mathscr{G}_{7,-2}$; see Fig. 6.

Now we consider $v_4\notin N_G(v_1)$, thus, $d_G(z_1)=a+b-1$ with
$z_1\not=v_4$. Applying (\ref{eq:1.1}) at $z_1$ yields
$d_G(z_2)=-ab-b^2+2b-1\in\{2,3,a+b-1\}$. It is routine to check that
only $d_G(z_2)=3$ holds, which implies $b=-1$ since  $a+b-1>3$.
Therefore, $z_2=v_4$ and $d_G(v_4)=3$. Note that $(a,b)=(7,-1),
a+b-1=5$. Applying (\ref{eq:1.1}) at $v_2$, we get that
$d_G(u_1)+d_G(w_1)=8$. If $d_G(w_1)=5,d_G(u_1)=3$ with $u_1=v_3$,
then $d_G(w_2)=3$. Hence, $w_2=v_4$. Then applying (\ref{eq:1.1}) at
$v_4$ yields $d_G(s_{k_3-1})=1$, it's impossible. So we have
$d_G(u_1)=5,d_G(w_1)=3$ with $w_1=v_4$. Moreover, $d_G(u_2)=3$, so
$u_2=v_3$. Applying (\ref{eq:1.1}) at $v_3$ yields $v_4\in
N_G(v_3)$. Thus, we obtain the graph $G\cong
G_{39}\in\mathscr{G}_{7,-1}$; see Fig. 6.

\medskip\noindent
$\bullet$ $t=3$. Thus, $v_2,v_3,v_4\in N_G(v_1)$ and
$d_G(v_2)=d_G(v_3)=d_G(v_4)=3$. By (\ref{eq:3.44}) we have
$3a+b=18$. Then applying (\ref{eq:1.1}) at $v_2$ yields
$d_G(u_1)+d_G(w_1)=3a+b-12=6$. Note that
$d_G(u_1),d_G(w_1)\in\{2,3,a+b-1\}$, it is routine to check that
$d_G(u_1),d_G(w_1)\not=2$. Together with $d_G(u_1)+d_G(w_1)=6$, we
have $d_G(u_1)=d_G(w_1)=3.$

First we assume $u_1=v_3,w_1=v_4$. Applying (\ref{eq:1.1}) at $v_3$,
we get $d_G(s_1)=3$. If $s_1=v_4$, thus $G$ is a regular graph with
only one Q-main eigenvalue, a contradiction. So we have
$s_1\not=v_4$. Hence, $a+b-1=d_G(s_1)=3$. It implies that
$(a,b)=(7,-3)$. Applying (\ref{eq:1.1}) at $s_1$, we have
$d_G(s_2)=5\notin\{2,3\}$, a contradiction.

Now we assume, without loss of generality, that $w_1\not=v_4$. We
may also get a contradiction by a similar discussion as above.

\vspace{2mm} \textbf{Subcase 1.2}.\ $t_1=1$, that is $\{d_G(r_1),
d_G(r_2), d_G(r_3)\}=\{2, t\cdot 3, (2-t)\cdot (a+b-1)\},$ where
$t=0,1,2.$ \vspace{2mm}

In this case, we assume, without loss of generality, that
$d_G(x_1)=2$.

Applying (\ref{eq:1.1}) at $v_1$ yields
\begin{equation}\label{eq:71}
   2+3t+(a+b-1)(2-t)=3a+b-9.
\end{equation}
Applying (\ref{eq:1.1}) at $x_1$, we get that
$d_G(x_2)=2a+b-7\in\{2,3,a+b-1\}$.

It's routine to check that only $t=2$ with $d_G(x_2)=a+b-1$ holds.
It implies that $d_G(v_3)=d_G(v_4)=3$ and $(a,b)=(6,-1)$. Hence,
$d_G(x_2)=4$.

First we consider $x_2\not=v_2$. Applying (\ref{eq:1.1}) at $x_2$
yields $d_G(x_3)=3$. Hence, we have $x_3=v_2$ with $d_G(v_2)=3$.  We
apply (\ref{eq:1.1}) at $v_2$ to get $d_G(u_1)+d_G(w_1)=4$. It
implies that $d_G(u_1)=d_G(w_1)=2$. Applying (\ref{eq:1.1}) at
$u_1$, we have $d_G(u_2)=4$, moreover, $d_G(u_3)=3$ with $u_3=v_3$.
Then applying (\ref{eq:1.1}) at $v_3$ yields $d_G(s_1)=1$, it's
impossible.

Now we consider $x_2=v_2$. Thus, $d_G(v_2)=4$. Applying Lemma
1.8(ii) at $v_2$ yields $d_G(u_1)+d_G(w_1)=4$, which implies that
$d_G(u_1)=d_G(w_1)=2$. Applying (\ref{eq:1.1}) at $u_1$ and $w_1$
respectively, we have $d_G(u_2)=d_G(w_2)=3$ with $u_2=v_3$ and
$w_2=v_4$. Then applying (\ref{eq:1.1}) at $v_3$ gives $v_4\in
N_G(v_3)$. Thus, it's simple to check that $G\cong
G_{40}\in\mathscr{G}_{6,-1}$; see Fig. 6.

\vspace{2mm} \textbf{Subcase 1.3}. $t_1=2$, that is $\{d_G(r_1),
d_G(r_2), d_G(r_3)\}=\{2,2,3\}$, or $\{2,2,a+b-1\}$. \vspace{2mm}

In this case, we assume, without loss of generality, that
$d_G(x_1)=d_G(y_1)=2$.

First consider $\{d_G(r_1), d_G(r_2), d_G(r_3)\}=\{2,2,3\}.$ Then
applying (\ref{eq:1.1}) at $v_1$ yields $2+2+3=3a+b-9$. It implies
that $3a+b=16$. Applying (\ref{eq:1.1}) at $x_1$ yields
$d_G(x_2)=2a+b-7\in\{2,3,a+b-1\}$. It's easy to check that
$d_G(x_2)\not=2$. For $d_G(x_2)\in\{3,a+b-1\}$, we have
$(a,b)=(6,-2)$ and $a+b-1=3$. If $x_2\not=v_2$, then applying
(\ref{eq:1.1}) at $x_2$ yields $d_G(x_3)=4\notin\{2,3\}$, a
contradiction. Hence, $x_2=v_2$. By a similar proof, we finally
obtain  $G\cong G_{24}\in\mathscr{G}_{6,-2}$, which has no pendant
vertices. It's a contradiction.

Now we consider $\{d_G(r_1), d_G(r_2), d_G(r_3)\}=\{2,2,a+b-1\}$.
Then applying (\ref{eq:1.1}) at $v_1$ yields $2+2+a+b-1=3a+b-9$,
which gives $a=6$. Since $a+b-1\geq3$, we have $b=-1,$ or $-2$.

If $(a,b)=(6,-2)$, then $a+b-1=3$. Thus, $\{d_G(r_1), d_G(r_2),
d_G(r_3)\}=\{2,2,3\}.$ By a similar discussion as above, this is
impossible.

If $(a,b)=(6,-1)$, then $d_G(z_1)=a+b-1=4$. If $z_1=v_4$, then
applying (\ref{eq:1.1}) at $v_4$ yields
$d_G(s_{k_3-1})+d_G(w_{k_2-1})=3$, it's impossible. Hence,
$z_1\not=v_4$. By (\ref{eq:1.1}) we can get that
$d_G(z_2)=2,d_G(z_3)=3$. Thus, $z_3=v_4$ with $d_G(v_4)=3$. Applying
(\ref{eq:1.1}) at $v_4$ yields $d_G(s_{k_3-1})+d_G(w_{k_2-1})=6$.
Hence, either $d_G(s_{k_3-1})=d_G(w_{k_2-1})=3$, or
$\{d_G(s_{k_3-1}), d_G(w_{k_2-1})\}=\{2, 4\}.$

If $d_G(s_{k_3-1})=d_G(w_{k_2-1})=3$, then $v_2,v_3\in N_G(v_4)$ and
$d_G(v_2)=d_G(v_3)=3$. Applying (\ref{eq:1.1}) at $x_1$ yields
$d_G(x_2)=4$, moreover, $d_G(x_3)=3$. Hence, $x_3=v_2$. Similarly,
$d_G(y_2)=4,d_G(y_3)=3$, and $y_3=v_3$. Applying (\ref{eq:1.1}) at
$v_2$ yields $d_G(u_1)=1$, a contradiction.

Hence, we consider $\{d_G(s_{k_3-1}), d_G(w_{k_2-1})\}=\{2, 4\}.$
Without loss of generality, we assume
$d_G(s_{k_3-1})=2,d_G(w_{k_2-1})=4$. Then applying (\ref{eq:1.1}) at
$w_{k_2-1}$ yields $d_G(w_{k_2-2})=2$, moreover, $d_G(w_{k_2-3})=3$.
Thus, we have $w_{k_2-3}=v_2$ with $d_G(v_2)=3$. We apply
(\ref{eq:1.1}) at $x_1$ to get $d_G(x_2)=4,$ moreover, $d_G(x_3)=3$
with $x_3=v_2$. Then applying (\ref{eq:1.1}) at $v_2$ yields
$d_G(u_1)=2$. Note that $d_G(u_1)=d_G(y_1)=d_G(s_{k_3-1})=2$, then
applying (\ref{eq:1.1}) at  $u_1,\,y_1,\,s_{k_3-1}$ respectively
yields $d_G(u_2)=d_G(y_2)=d_G(s_{k_3-2})=4$.

If $u_2\not=v_3$, then applying (\ref{eq:1.1}) at $u_2$ yields
$d_G(u_3)=3$, which implies that $u_3=v_3$ and $d_G(v_3)=3$. Thus,
we have $y_2\not=v_3,\,s_{k_3-2}\not=v_3$. By a similar discussion
as $u_2$, we may get that $y_3=s_{k_3-3}=v_3$. However, by applying
(\ref{eq:1.1}) at $v_3$, we have $4+4+4=3a+b-9$ with $(a,b)=(6,-1)$.
It's a contradiction. Therefore, $u_2=v_3$. Similarly,
$y_2=s_{k_3-2}=v_3$. Thus, it's easy to check that $G\cong
G_{41}\in\mathscr{G}_{6,-1}$; see Fig. 6.

\vspace{2mm} \textbf{Subcase 1.4}. $t_1=3$, that is
$d_G(r_1)=d_G(r_2)= d_G(r_3)=2$. \vspace{2mm}

In this case, $d_G(x_1)=d_G(y_1)=d_G(z_1)=2$.

Applying (\ref{eq:1.1}) at $v_1$, we have $2+2+2=3a+b-9$, i.e.,
$3a+b=15$. Applying (\ref{eq:1.1}) at $x_1$ yields
$d_G(x_2)=2a+b-7\in\{2,3,a+b-1\}$. It's easy to check that this is
impossible.

\vspace{2mm} {\bf Case 2.}\ Each of the vertices in
$\{v_1,\,v_2,\,v_3,\,v_4\}$ is an attached-vertex. That is,
$d_G(v_1)=d_G(v_2)=d_G(v_3)=d_G(v_4)=a+b-1>3$. Let
$N_{\widetilde{G}}(v_1)=\{r_1, r_2, r_3\}$, then $\{d_G(r_1),
d_G(r_2), d_G(r_3)\}=\{t\cdot 2, (3-t)\cdot (a+b-1)\},$ where $0\leq
t\leq3$.

Applying Lemma \ref{lem4}(ii) at $v_1$, we have
$d_G(r_1)+d_G(r_2)+d_G(r_3)=-ab-b^2+2b+3$. Hence,
\begin{equation}\label{eq:70}
3t+(a+b-1)(3-t)=-ab-b^2+2b+3.
\end{equation}

\medskip\noindent
$\bullet$ $t=0$. By (\ref{eq:70}) we have $ab+b^2+3a+b=6$. If
$v_2,v_3,v_4\in N_G(v_1)$, we apply (\ref{eq:1.1}) at $v_2$ to get
$d_G(u_1)+d_G(w_1)=-ab-b^2-a+b+4$. Note that
$d_G(u_1),d_G(w_1)\in\{2,a+b-1\}$. It is routine to check that only
$d_G(u_1)=d_G(w_1)=a+b-1$ is true.

First we assume that $u_1=v_3,w_1=v_4$. Applying (\ref{eq:1.1}) at
$v_3$, we have $d_G(s_1)=a+b-1$. By (\ref{eq:1.1}), it's easy to
check that $s_1=v_4$. Hence, we obtain the graph $G\cong
G_{42}\in\mathscr{G}(a,b)$, where $a, b$ satisfying $ab+b^2+3a+b=6$,
with $a+b-1\ge 4$ and $b\le -1$. Denote the number of pendant
vertices attached at $v_i(i=1,2,3,4)$ by k, then
%In particular, $a+b-1=\frac{6-3a-2b}{b}=\frac{6-3a}{b}-2\geq4$, which gives that
$k=a+b-1-3\geq1$; see Fig. 6.

Now we assume, without loss of generality, $w_1\neq v_4$. Applying
(\ref{eq:1.1}) at $w_1$, we get that
$d_G(w_2)=-ab-b^2-a+b+3\in\{2,a+b-1\}$, which gives no integer
solution such that $a+b-1>3$, contradiction.

If $\{v_2,v_3,v_4\}$ contains a member, say $v_4$, not in
$N_G(v_1)$. Thus, $d_G(z_1)=a+b-1$ with $z_1\not=v_4$. Applying
(\ref{eq:1.1}) at $z_1$ yields
$d_G(z_2)=-ab-b^2-a+b+3\in\{2,a+b-1\}$. It is routine to check that
this is impossible.

\medskip\noindent
$\bullet$ $t=1,2,3$. In this subcase, we assume,  without loss of
generality, that $d_G(x_1)=2$. Applying (\ref{eq:1.1}) at $x_1$, we
get that $d_G(x_2)=a-3\in\{2,a+b-1\}$. It's easy to check that
$d_G(x_2)\not=2$. So we consider $d_G(x_2)=a+b-1$. It implies that
$b=-2$. For $t=1$, (\ref{eq:70}) gives $ab+b^2+2a=3$; for $t=2$,
(\ref{eq:70}) gives $ab+b^2+a-b=0$; for $t=3$, (\ref{eq:70}) gives
$ab+b^2-2b+3=0$. Each gives no integer solution such that $a+b-1>3$
since $b=-2$. It's a contradiction.

Thus, we complete the proof.
\end{proof}
\begin{thm}
$G_{28}\in \mathscr{G}_{6,-1}, G_{29}\in \mathscr{G}_{6,-1}, G_{30}\in\mathscr{G}_{6,-1}, G_{31}\in\mathscr{G}_{6,-1}, G_{32}\in\mathscr{G}_{7,-2}, G_{33}\in\mathscr{G}_{6,-1},
G_{34}\in \mathscr{G}_{6,-2}, G_{35}\in\mathscr{G}_{6,-2}, G_{36}\in\mathscr{G}_{7,-1}, G_{37}\in\mathscr{G}_{8,-2}, G_{38}\in\mathscr{G}_{7,-2},
G_{39}\in\mathscr{G}_{7,-1}, G_{40}\in\mathscr{G}_{6,-1}, G_{41}\in\mathscr{G}_{6,-1}, G_{42}\in\mathscr{G}_{a,b}$ (see Figs. 3-6) are all the tricyclic graphs with pendants having exactly two
$Q$-main eigenvalues.
\end{thm}
\begin{proof}
By Propositions 2,3 and 5, $G\not\in \mathscr{G}_{a,b}$ if $\widetilde{G}=T_3, T_9, T_{13}$. In view of the proof of Proposition 1, we obtain that 
$G\not\in \mathscr{G}_{a,b}$ if $\widetilde{G}=T_1, T_2, T_5, T_7, T_8, T_{10}$. 
Hence, our results follow immediately from Propositions 1, 4, 6, and 7.
\end{proof}

\end{document}